\documentclass[11pt]{article}

\usepackage[utf8]{inputenc}
\usepackage{geometry}
\usepackage{amsmath,amssymb,amsfonts,mathtools}
\usepackage{mathrsfs}
\usepackage{graphicx}
\usepackage{subcaption}
\usepackage{float}
\usepackage[section]{placeins}

\usepackage{booktabs}
\usepackage{array}
\usepackage{multirow}
\usepackage{siunitx}

\usepackage{algorithm}
\usepackage{algpseudocode}

\usepackage{appendix}
\usepackage{chngcntr} 

\usepackage{enumitem}
\usepackage[table]{xcolor}

\usepackage[bookmarks=false]{hyperref}
\usepackage{url}
\usepackage{cite}

\usepackage{authblk}
\usepackage{bm}
\newcommand{\bx}{\mathbf{x}}

\newcommand{\incfig}[2][]{%
  \IfFileExists{#2}{\includegraphics[#1]{#2}}{\fbox{\texttt{Missing: \detokenize{#2}}}}%
}

\usepackage[table]{xcolor}

\colorlet{highlightblue}{blue}   % 蓝色修改版
\title{Random Projection Neural Networks for Phase-Field Models: Model Reformulation, Scaling, and Network Architecture}
\author[1]{Chunyan Li}

\author[1]{Qi Wang\thanks{Corresponding author: \texttt{qwang@math.sc.edu}}}
\affil[1]{Department of Mathematics, University of South Carolina, Columbia, SC 29208, USA}
\date{}

\begin{document}
\maketitle

\begin{abstract}
Phase-field models, particularly the Allen–Cahn (AC) and Cahn–Hilliard (CH) equations, are widely used to describe interfacial dynamics, but resolving their sharp transition layers and stiff dynamics presents substantial numerical challenges. Physics-informed neural networks (PINNs) offer a mesh-free approach, yet optimizing large numbers of network parameters can be computationally demanding and may still yield limited accuracy. We develop a random projection neural network (RPNN) framework that combines first-order PDE reformulations with independent random-feature representations of the unknowns. By fixing randomly initialized hidden-layer parameters and optimizing only the output coefficients, RPNNs reduce training to  least-squares problems amenable to Newton-type solvers. Introducing auxiliary variables eliminates higher-order derivatives from the PDE residuals, while separate random feature bases in the unknowns provide flexibility in approximating fields with different scales and regularity. For the CH equation, we additionally scale the chemical potential to address the disparity between its magnitude and that of the phase field. We compare the first-order and corresponding higher-order formulations and systematically examine the effects of initialization distributions, the random feature range \(R_\xi\), seven activation functions, and dense versus separable space-time bases. We also evaluate residual-based adaptive refinement, energy-dissipation-rate-guided and gradient-based adaptive sampling, residual reweighting, and mild spatial derivative-penalty regularization. Across the reported benchmarks, the proposed framework reduces relative \(L^2\) and absolute \(L^\infty\) errors by one to three orders of magnitude compared with the tested PINN baselines under matched conditions. These results demonstrate the potential of first-order RPNNs for accurate approximation of phase-field dynamics with sharp transition layers.
\end{abstract}


 %  For the CH equation, it achieves a relative \(L^2\) error of \(8.64\times10^{-4}\) and an absolute \(L^\infty\) error of \(1.39\times10^{-3}\), while providing a continuous spatio-temporal representation. 
% \input{introduction}
\section{Introduction}\label{sec:introduction}

Phase-field models, particularly the Allen--Cahn (AC) and Cahn--Hilliard (CH) equations, are widely used to describe phase separation, coarsening dynamics, and interfacial motion~\cite{chen2002phase,Zhang&W2010,Li&W2014}. Their diffuse-interface descriptions avoid explicit interface tracking, but accurately resolving thin transition layers and stiff dynamics remains computationally challenging, particularly over long time intervals or in systems with dynamic boundary conditions in complex domains~\cite{jing2022entropy,Jing2023cms}. Traditional finite difference, finite element, and spectral methods provide well-established approaches to these problems. Nevertheless, thin interfaces can require fine spatial resolution, while stiffness can impose restrictive time steps on explicit schemes or require the solution of coupled nonlinear systems in implicit methods. These challenges motivate the development of alternative approximation strategies for phase-field dynamics.

Physics-informed neural networks (PINNs)~\cite{raissi2019physics} provide a mesh-free approach by representing PDE solutions with neural networks and incorporating governing equations, initial conditions, and boundary conditions into a residual-based loss evaluated at collocation points. Automatic differentiation supplies the derivatives required by the PDE residuals, while the neural representation offers flexibility on irregular domains. Several strategies have been developed to improve their performance for phase-field problems. Temporal approaches include adaptive time stepping and sequential training over time segments~\cite{wight2020solving,mattey2022novel}, time-iteration schemes~\cite{li2022delisa}, and discrete-time formulations with adaptive collocation or Runge--Kutta integration~\cite{xu2022adaptive,he2023adaptive}. Other developments include pseudo-spectral PINNs for identifying phase-field models from image data~\cite{zhao2021discovering}, multi-output neural representations~\cite{singh2023multioutput}, and fully discrete convolutional neural network solvers~\cite{geng2024deep}. Adaptive loss-weighting methods offer another route to improving residual-based training~\cite{mcclenny2023self,liu2024discontinuity,zhang2026priori}. Physical structure has also been incorporated through energy-based sequential training~\cite{guo2024adaptive}, energy-dissipation-rate-guided adaptive sampling (EDRAS)~\cite{li2025edras}, and energy-dissipation-preserving loss terms~\cite{kutuk2025energy}. Despite these advances, achieving high accuracy with PINNs remains challenging for phase-field problems with sharp transition layers, especially, when high spatial derivatives are present in the models such as in the Cahn-Hilliard model. Training large numbers of network parameters involves highly nonconvex optimization and can be sensitive to network architecture, initialization, loss weighting, and sampling strategies. These difficulties motivate approaches that reduce the number of trainable parameters and simplify training.

Randomized neural architectures, whose origins date back to 1992~\cite{schmidt1992elm,huang2004extreme,huang2011extreme,huang2015extreme}, provide an alternative to fully trainable deep neural networks. Their hidden-layer weights and biases are typically sampled from a prescribed distribution, such as $\mathcal{U}(-R_\xi,R_\xi)$ with random feature range $R_\xi$, and then fixed, leaving only the output-layer coefficients to be determined. For residual-based PDE solvers, this construction recasts neural network training as a least-squares problem in those coefficients. In particular, it becomes a linear programming for linear PDE problems. The resulting problem can be solved using direct least-squares methods or Newton-type algorithms, including Gauss--Newton and Levenberg--Marquardt, which exploit approximate second-order information rather than relying on the first-order gradient-based optimizers commonly used to train PINNs. Closely related approaches appear under different names, including extreme learning machines (ELMs), random feature methods (RFMs), and random projection neural networks (RPNNs). Representative developments include physics-informed ELM (PIELM) for linear PDEs~\cite{dwivedi2020pielm}; rank-inspired neural networks (RINNs), which augment direct random initialization with covariance-based preconditioning of the resulting hidden-layer feature matrix~\cite{peng2026RINN}; RFMs for stationary and time-dependent PDEs and interface problems~\cite{chen2022bridging,chen2023random,chi2024random,liao2025solving,mei2025solving,linghu2026liquid}; ELM-based methods for time-dependent problems~\cite{dong2021local,dong2021modified,dong2022computing,wang2024extreme,calabro2021extreme,he2024new}; and RPNNs for function approximation and the numerical solution of ODEs and PDEs~\cite{fabiani2023parsimonious,fabiani2025random,fabiani2026stability,georgiou2025heatnets,van2026local}. Theoretical studies have also established universal approximation properties and, under suitable assumptions, error rates comparable to those of Legendre polynomial approximation~\cite{fabiani2025random}, together with extensions to operator learning~\cite{fabiani2025randonets}.

Random-feature and ELM-based methods have already been applied to phase-field problems~\cite{he2024new,liao2025solving,linghu2026liquid,zhou2026discrete}. Their attainable accuracy, however, depends on both the PDE formulation and the approximation spaces generated by the random features. This interaction is particularly important for the fourth-order CH equation, where higher-order derivatives and coupled fields with different magnitudes and spatial structures can place competing demands on a single feature representation. At the same time, feature quality depends on the activation function, initialization distribution, random feature range $R_\xi$, and space-time basis construction. These considerations raise two related questions: \textit{how should the phase-field equations and their coupled fields be represented to obtain accurate RPNN approximations, and which random-feature choices most strongly affect the resulting performance}? Addressing these questions requires a systematic assessment of equation formulation, field representation, and feature design within a common numerical framework.

In this work, we develop an RPNN least-squares collocation framework based on first-order formulations of the AC and CH equations. This choice is motivated by the variationally correct neural residual regression framework of Bachmayr, Dahmen, and Oster~\cite{bachmayr2025variational}, which establishes, for suitable linear PDE problems, stable first-order formulations with computable residual losses equivalent to the squared solution error in appropriate norms. This equivalence makes the loss a reliable measure of approximation error rather than merely an optimization objective: a small loss guarantees a small solution error, and the current loss provides a posteriori accuracy control at any stage of training, independently of whether optimization has reached a minimizer. These results provide a principled motivation for investigating first-order residual formulations for nonlinear phase-field equations, although establishing the corresponding error guarantees requires separate analysis and is left for a future investigation. Introducing auxiliary variables also replaces higher-order derivatives of the order parameter with first derivatives of independently approximated fields. Within an RPNN framework, this permits separate random feature spaces and field-specific approximation budgets. We investigate this combination numerically for nonlinear phase-field dynamics, with particular attention to the scaling of the chemical potential in the CH system. % Establishing residual--error norm equivalence for these nonlinear AC and CH formulations is beyond the scope of the present work and is left for future investigation.

Our contributions in this study are fourfold. First, we construct  RPNN approximations using independent random feature bases for the solution fields in first-order PDEs and compare them with the corresponding higher-order formulations. For the CH equation, we examine three equivalent auxiliary-variable formulations and demonstrate that appropriate chemical-potential scaling improves outcomes of numerical optimization. Independent bases also provide greater flexibility than a shared one, allowing the chemical potential to be approximated with substantially fewer basis functions while maintaining or improving accuracy. Second, we systematically assess the influence of initialization distributions, the random feature range $R_\xi$, seven activation functions, and dense versus separable space-time bases, identifying their effects on approximation accuracy. Third, we propose gradient-based adaptive sampling, residual reweighting, and mild derivative-penalty regularization, and compare them with other adaptive refinement methods—including residual-based adaptive refinement (RAR)~\cite{lu2021deepxde}, EDRAS~\cite{li2025edras}, and their combination—to determine their benefits and limitations for resolving sharp transition layers in phase filed models. Finally, numerical experiments on one-dimensional AC and CH benchmarks demonstrate one- to three-order-of-magnitude improvements in relative $L^2$ and absolute $L^\infty$ errors compared with the tested PINN baselines and five state-of-the-art methods in the literature. 

The remainder of this paper is organized as follows. Section~\ref{sec:method} presents the numerical framework, including the phase-field models and their reformulations in Section~\ref{subsec:models}, standard RPNN representations in Section~\ref{subsec:RPNN}, domain decomposition in Section~\ref{subsec:local_ddm}, nonlinear least-squares training in Section~\ref{subsec:least_squares}, and adaptive sampling and residual reweighting in Section~\ref{subsec:adap_samp}. Sections~\ref{sec:AC1d} and~\ref{sec:results_ch} report numerical experiments for the AC and CH equations, respectively, examining the effects of equation formulation, scaling, random-feature design, and enhancement strategies. Section~\ref{sec:conclusion} summarizes the findings, discusses limitations, and outlines future research directions. 
\section{Methodology}\label{sec:method}

\subsection{Phase-field models and first-order formulations}\label{subsec:models}
Let $Q=\Omega\times(0,T]$, where $\Omega\subset\mathbb R^d$ is the spatial domain. We consider the Ginzburg--Landau free energy
\begin{equation}\label{eq:GL_energy}
E[u]=\int_\Omega\left[\frac{\gamma_1}{2}|\nabla u|^2
+\frac{\gamma_2}{4}(u^2-1)^2\right]\,\mathrm d\mathbf x,
\qquad \gamma_1,\gamma_2>0.
\end{equation}
Here $\gamma_1$ penalizes spatial variations, $\gamma_2$ determines the double-well potential, and $\ell=\sqrt{\gamma_1/\gamma_2}$ characterizes the diffuse-interface thickness. Under homogeneous Neumann boundary conditions, the variational derivative of the energy is the chemical potential
\begin{equation}\label{eq:mu}
\mu=\frac{\delta E}{\delta u}
=-\gamma_1\Delta u+\gamma_2(u^3-u).
\end{equation}
% For sufficiently regular solutions, $\mathrm dE/\mathrm dt=\int_\Omega\mu u_t\,\mathrm d\mathbf x$.

\paragraph{Allen--Cahn equation.}
With unit mobility, the Allen--Cahn equation is the $L^2$ gradient flow of \eqref{eq:GL_energy}:
\begin{equation}\label{eq:AC}
u_t=-\mu=\gamma_1\Delta u-\gamma_2(u^3-u)
\qquad\text{in }Q,
\end{equation}
with $u(\mathbf x,0)=u_0(\mathbf x)$ and $\nabla u\cdot\mathbf n=0$ on $\partial\Omega\times(0,T]$. Introduce an independent auxiliary field $\mathbf q=\nabla u$. The equivalent first-order system is
\begin{equation}\label{eq:AC_first_order}
\begin{aligned}
u_t-\gamma_1\nabla\cdot\mathbf q+\gamma_2(u^3-u)&=0,\\
\mathbf q-\nabla u&=0,
\end{aligned}
\end{equation}
subject to
\begin{equation}\label{eq:AC_first_order_data}
u(\mathbf x,0)=u_0(\mathbf x),\qquad
\mathbf q\cdot\mathbf n=0\quad\text{on }\partial\Omega\times(0,T].
\end{equation}
The energy dissipation law is
\begin{equation}\label{eq:AC_energy_dissipation}
\frac{\mathrm dE}{\mathrm dt}
=-\int_\Omega\mu^2\,\mathrm d\mathbf x
=-\int_\Omega u_t^2\,\mathrm d\mathbf x\leq0.
\end{equation}

\paragraph{Cahn--Hilliard equation.}
The Cahn--Hilliard equation is the mass-conserving $H^{-1}$ gradient flow of the same energy. With unit mobility, its mixed formulation is
\begin{equation}\label{eq:CH}
\begin{aligned}
u_t-\Delta\mu&=0,\\
\mu+\gamma_1\Delta u-\gamma_2(u^3-u)&=0,
\end{aligned}
\qquad\text{in }Q,
\end{equation}
with $u(\mathbf x,0)=u_0(\mathbf x)$ and $\nabla u\cdot\mathbf n=\nabla\mu\cdot\mathbf n=0$ on the spatial boundary. Introducing $\mathbf q=\nabla u$ and $\mathbf j=\nabla\mu$ yields
\begin{equation}\label{eq:CH_first_order}
\begin{aligned}
u_t-\nabla\cdot\mathbf j&=0,\\
\mu+\gamma_1\nabla\cdot\mathbf q-\gamma_2(u^3-u)&=0,\\
\mathbf q-\nabla u&=0,\\
\mathbf j-\nabla\mu&=0,
\end{aligned}
\end{equation}
with
\begin{equation}\label{eq:CH_first_order_data}
u(\mathbf x,0)=u_0(\mathbf x),\qquad
\mathbf q\cdot\mathbf n=\mathbf j\cdot\mathbf n=0
\quad\text{on }\partial\Omega\times(0,T].
\end{equation}
The physical mass flux is $-\mathbf j$. The exact solution conserves $\int_\Omega u\,\mathrm d\mathbf x$ and satisfies
\begin{equation}\label{eq:CH_energy_dissipation}
\frac{\mathrm dE}{\mathrm dt}
=-\int_\Omega|\nabla\mu|^2\,\mathrm d\mathbf x\leq0.
\end{equation}

\paragraph{Chemical-potential scaling and one-dimensional implementation.}
To adjust the relative magnitudes of the approximated fields, let $\mu=s_\mu\widetilde\mu$ and $\mathbf j=s_\mu\widetilde{\mathbf j}$, with a prescribed constant $s_\mu>0$. The scaled first-order CH system becomes
\begin{equation}\label{eq:CH_first_order_scaled}
\begin{aligned}
u_t-s_\mu\nabla\cdot\widetilde{\mathbf j}&=0,\\
\widetilde\mu+\frac{\gamma_1}{s_\mu}\nabla\cdot\mathbf q
-\frac{\gamma_2}{s_\mu}(u^3-u)&=0,\\
\mathbf q-\nabla u&=0,\\
\widetilde{\mathbf j}-\nabla\widetilde\mu&=0.
\end{aligned}
\end{equation}
For example, $s_\mu=\gamma_2$ normalizes the coefficient of the nonlinear term and gives $\gamma_1/s_\mu=\ell^2$. The initial condition on $u$ is unchanged, and the boundary conditions are $\mathbf q\cdot\mathbf n=\widetilde{\mathbf j}\cdot\mathbf n=0$. Although scaling preserves the zero-residual PDE system, it changes the relative weighting of residual blocks unless their loss weights are transformed accordingly. The CH formulations and scaling choices used in the numerical comparisons are specified in Section~\ref{sec:results_ch}.

%All numerical experiments considered here are one-dimensional. Thus $\nabla=\partial_x$, $\nabla\cdot=\partial_x$, and the auxiliary fields reduce to scalars $q=u_x$ and $j=\mu_x$. 
The AC first-order system has two scalar unknowns $(u,q)$, whereas the CH system has four, $(u,q,\mu,j)$, or their scaled counterparts. Only $u$ requires independent initial data; auxiliary initial values follow from the defining relations for sufficiently smooth solutions. We approximate each field independently and penalize all first-order equation residuals. Consequently, evaluating these residuals requires only first derivatives of the network outputs. The second-order AC equation and the higher-order CH formulations are retained for numerical comparison.

% --------------------------------------------------
\subsection{Random projection neural networks}\label{subsec:RPNN}

\paragraph{Input normalization.} 
On a local space-time domain $[a_k,b_k]\times[t_n,t_{n+1}]$, normalize the physical coordinates by
\begin{equation}\label{eq:input_normalization}
\widehat x=\frac{2(x-a_k)}{b_k-a_k}-1,
\qquad
\widehat t=\frac{2(t-t_n)}{t_{n+1}-t_n}-1.
\end{equation}
In higher dimensions, the spatial normalization is applied componentwise on rectangular subdomains. This coordinate transformation introduces no trainable parameters. Physical derivatives include the chain-rule factors
\begin{equation}\label{eq:normalization_derivatives}
\partial_x=\frac{2}{b_k-a_k}\partial_{\widehat x},
\qquad
\partial_t=\frac{2}{t_{n+1}-t_n}\partial_{\widehat t}.
\end{equation}

\paragraph{Random projection neural network.} 
Each scalar field $f$ in the chosen PDE formulation is approximated by a single-hidden-layer RPNN:
\begin{equation}\label{eq:rpnn_general}
f(\mathbf x,t;\mathbf w_f)
=\sum_{j=1}^{m_f}w_{f,j}H_{f,j}
(\widehat{\mathbf x},\widehat t),
\end{equation}
where $H_{f,j}$ denotes a random feature function whose hidden-layer weights and biases are sampled once and then held fixed. Only the output coefficients $\mathbf w_f$ are optimized. In the independent-basis construction, each field uses its own random feature realization and may have a different number of features $m_f$. In the shared-basis comparison, the fields use the same feature functions but retain separate output coefficients. The feature constructions and the distributions and scales used to sample their hidden parameters are specified below.

\paragraph{Dense and separable space-time features.}
We consider dense features, which combine the normalized spatial and temporal coordinates through a single affine transformation followed by an activation function~\cite{chen2022bridging,chen2026physics}:
\begin{equation}\label{eq:dense_features}
H_{f,j}(\widehat{\mathbf x},\widehat t)
=\sigma\!\left(\mathbf a_{f,j}\cdot\widehat{\mathbf x}
+b_{f,j}\widehat t+c_{f,j}\right).
\end{equation}
This construction couples space and time within each feature, while allowing their parameters to be sampled using different distributions or scales.

We compare this construction with separable features motivated by separation of variables~\cite{chen2023random,chen2026physics}:
\begin{equation}\label{eq:sep_features}
H_{f,j}(\widehat{\mathbf x},\widehat t)
=\sigma_x\!\left(\mathbf a_{f,j}\cdot\widehat{\mathbf x}+c_{f,j}\right)
\sigma_t\!\left(b_{f,j}\widehat t+d_{f,j}\right).
\end{equation}
Each separable feature is a product of spatial and temporal factors, whose activation functions may coincide or differ. For both constructions, the required feature derivatives are evaluated using the chain rule, enabling direct assembly of the PDE residuals and their Jacobian with respect to the output coefficients.

\paragraph{Initialization distributions and feature scales.}
We compare uniform and normal initialization. For a scalar hidden parameter $\xi$ with assigned scale $R_\xi>0$, use
\begin{equation}\label{eq:feature_initialization}
\xi\sim\mathcal U(-R_\xi,R_\xi)
\qquad\text{or}\qquad
\xi\sim\mathcal N\!\left(0,\frac{R_\xi^2}{3}\right).
\end{equation}
All hidden-layer weights and biases are sampled independently. The spatial weights, temporal weights, and biases may have distinct sampling scales, denoted by $R_\xi^{(\mathbf x)}$, $R_\xi^{(t)}$, and $R_\xi^{(c)}$, respectively, with an additional bias scale $R_\xi^{(d)}$ for separable features. These scales may also vary between fields. The distributions in \eqref{eq:feature_initialization} have matching variances, allowing uniform and Gaussian initialization to be compared at the same second-moment scale. For Gaussian initialization, $R_\xi$ specifies the scale rather than a bound on sampled values.

The feature scale controls spatial or temporal variation in an activation-dependent manner: larger slopes produce more rapid oscillations for oscillatory activations or steeper transitions for sigmoidal activations. Since the appropriate scale depends on the PDE, activation, and approximation space~\cite{dong2022computing,fabiani2025random}, initialization distributions, $R_\xi$, activation functions, and dense versus separable features are assessed systematically in the numerical experiments.

% --------------------------------------------------
\subsection{Local models and domain decomposition}\label{subsec:local_ddm}
We combine sequential time-block solving with non-overlapping space--time domain decomposition following locELM~\cite{dong2021local}. First, divide the full time interval into $N_b$ blocks,
\[
0=T_0<T_1<\cdots<T_{N_b}=T,
\]
which are solved sequentially~\cite{wight2020solving,mattey2022novel,li2025edras}.
Within each block $\Omega\times[T_m,T_{m+1}]$, partition space and time as
\[
a=x_0<\cdots<x_{N_x}=b,
\qquad
T_m=t_0^{(m)}<\cdots<t_{N_t}^{(m)}=T_{m+1},
\]
giving local subdomains
\[
Q_{kn}^{(m)}
=(x_k,x_{k+1})\times(t_n^{(m)},t_{n+1}^{(m)}].
\]
Each local field is represented by \eqref{eq:rpnn_general}, with local coordinate normalization \eqref{eq:input_normalization} and subdomain-specific output coefficients. Within each time block, all local coefficients are solved jointly through a coupled least-squares
problem. PDE residuals are evaluated within each subdomain, and physical boundary conditions are imposed on the external spatial boundary.

For the first-order formulations, let $\mathbf U_{kn}^{(m)}$ collect the local fields: $(u,q)$ for AC and $(u,q,\mu,j)$ for CH, with scaled variables substituted when applicable. Within each time block, all field values are matched across adjacent spatial and temporal subdomains:
\begin{align}
\mathbf U_{kn}^{(m)}(x_{k+1},t)
-\mathbf U_{k+1,n}^{(m)}(x_{k+1},t)
&=\mathbf0,
\label{eq:c0x}\\
\mathbf U_{k,n+1}^{(m)}(x,t_{n+1}^{(m)})
-\mathbf U_{kn}^{(m)}(x,t_{n+1}^{(m)})
&=\mathbf0.
\label{eq:c0t_internal}
\end{align}
No additional derivative matching is imposed for these first-order systems. For the second-order AC and mixed second-order CH comparisons, field values are likewise matched across spatial and temporal interfaces, with additional continuity of the first spatial derivatives:
\begin{equation}\label{eq:c1x}
\partial_x f_{kn}^{(m)}(x_{k+1},t)
-\partial_x f_{k+1,n}^{(m)}(x_{k+1},t)=0,
\end{equation}
where $f=u$ for AC and $f\in\{u,\mu\}$ for mixed CH.

Between sequential time blocks, only the terminal approximation of $u$ is passed forward as initial data. Thus, the first block uses $u_0$, while subsequent blocks impose
\begin{equation}\label{eq:c0t}
u_{k,0}^{(m+1)}(x,T_{m+1})
-u_{k,N_t-1}^{(m)}(x,T_{m+1})=0,
\qquad m=0,\ldots,N_b-2.
\end{equation}
The terminal approximation from the preceding block is held fixed while the next block is solved. Auxiliary fields receive no independent initial data at block boundaries. All interface and initial matching conditions are included as least-squares residuals and are therefore satisfied approximately.

% Least-squares and enhancement descriptions are condensed below.
 
\subsection{Least-squares training problem}\label{subsec:least_squares}
Let $\mathbf w$ collect the output coefficients of all local fields.
We formulate the training problem as
\begin{equation}\label{eq:nlls_general}
\min_{\mathbf w}\mathcal L(\mathbf w)
=\frac12\|r(\mathbf w)\|_2^2,
\end{equation}
where the weighted residual vector is
\[
r(\mathbf w)=
\bigl[r_{\mathrm{PDE}}^\top,\,
\alpha_{\mathrm{IC}}r_{\mathrm{IC}}^\top,\,
\alpha_{\mathrm{BC}}r_{\mathrm{BC}}^\top,\,
\alpha_{\mathrm{IF}}r_{\mathrm{IF}}^\top\bigr]^\top.
\]
Here, $r_{\mathrm{PDE}}$ collects all equation residuals, including the auxiliary-variable relations, while $r_{\mathrm{IC}}$, $r_{\mathrm{BC}}$, and $r_{\mathrm{IF}}$ collect the initial, physical boundary, and spatial interface residuals, respectively. The factors $\alpha_{\mathrm{IC}}$, $\alpha_{\mathrm{BC}}$, and $\alpha_{\mathrm{IF}}$ control their relative contributions; their squares are the corresponding objective weights. For sequential time slabs, initial data are supplied by the preceding terminal approximation.

With the hidden-layer parameters fixed, the residual vector depends only on $\mathbf w$, and its Jacobian $J(\mathbf w)=\partial r/\partial\mathbf w$ can be assembled explicitly from the feature values and their derivatives. We solve \eqref{eq:nlls_general} using SciPy's \texttt{least\_squares} routine with this supplied Jacobian. To reduce sensitivity to the initial output coefficients, we incorporate the NLSQ-perturb strategy described in~\cite{dong2021local}, which supplements least-squares solves with perturbation-based restarts.
 
% ---------------------------------------------------------------
\subsection{Adaptive sampling and residual reweighting strategies}
\label{subsec:adap_samp}
We evaluate established adaptive sampling and residual reweighting methods, introduce a gradient-based refinement strategy (Grad-ADS), and compare two derivative-penalty regularization methods. Sampling-based refinement follows a two-stage procedure: the first solve provides an approximation from which refinement indicators are computed; the collocation set is then enriched, and the resulting coupled least-squares problem is solved using the first-stage output coefficients as the initial guess.

\paragraph{Residual-based and energy-dissipation-based refinement.}
Residual-based adaptive refinement (RAR) selects additional collocation points with the largest PDE residual indicator~\cite{lu2021deepxde}. Energy-dissipation-rate-based adaptive sampling (EDRAS) instead targets regions of high local energy
dissipation~\cite{li2025edras}, using $\rho=|u_t|^2$ for AC and $\rho=|\mu_x|^2$ for CH, where $\mu$ is the physical chemical potential. For scaled formulations, $\mu$ is recovered before evaluating the indicator. We also consider the combined score
\begin{equation}\label{eq:combined_indicator}
S_i=
\left(\frac{R_i}{R_{\mathrm{ref}}}\right)^{p_{\mathrm{res}}}
+\lambda_e
\left(\frac{\rho_i}{E_{\mathrm{ref}}}\right)^{p_e},
\end{equation}
where $R_i\geq0$ denotes the residual indicator. The normalization factors $R_{\mathrm{ref}}$ and $E_{\mathrm{ref}}$ are the root-mean-square values of the respective indicators over the candidate pool, plus $10^{-12}$. We use $p_{\mathrm{res}}=p_e=1$ and $\lambda_e=1$.

% ----------- residual reweighting -----------------------
\paragraph{Residual reweighting.}
We retain the collocation set and assign weights based on the first-stage PDE residuals. Let $R_i$ denote the residual at point $i$. We define
\begin{equation}\label{eq:residual_weights}
\omega_i=
1+\alpha\left(\frac{|R_i|}{R_{\mathrm{ref}}}\right)^p
\end{equation}
where
$R_{\mathrm{ref}}=\max\{\operatorname{median}_i|R_i|,10^{-12}\}$ is computed over all interior collocation points across subdomains. We use $\alpha=5$ and $p=2$. The second-stage PDE loss is
\[
\mathcal L_{\mathrm{PDE}}(\mathbf w)
=\frac12\sum_i\omega_i\,|r_i(\mathbf w)|^2,
\]
with $\omega_i$ held fixed during the solve. The initial, boundary, and interface terms remain unchanged.
% ------------- gradient based adaptive sampling --------------
\paragraph{Gradient-based adaptive sampling (Grad-ADS).}
We introduce Grad-ADS to target sharp transition layers using the spatial-gradient indicator
\[
G(x,t)=|\partial_xu(x,t)|.
\]
Points with indicator values exceeding a prescribed quantile threshold are selected, and their spatial locations are grouped into intervals. Each interval is expanded to form an interface band, which is then intersected with the local subdomains to determine the refinement regions. Additional collocation points are sampled uniformly within these regions and added to the existing sampling sets. This strategy identifies sharp interfaces even when they evolve slowly in time. Algorithm~\ref{alg:grad-refine} summarizes the procedure.

\paragraph{Derivative-penalty regularization.}
We also compare penalties on the first and second spatial derivatives:
\begin{equation}\label{eq:grad_penalty}
\mathcal L_{\mathrm{reg}}(\mathbf w)
=\mathcal L(\mathbf w)
+\frac{\lambda_D}{2}
\|D u(\cdot;\mathbf w)\|_{L^2(Q)}^2,
\qquad
D\in\{\partial_x,\partial_{xx}\}.
\end{equation}
Each discretized penalty is incorporated as an additional residual block with multiplier $\sqrt{\lambda_D}$, including the appropriate regularization weights. The penalty strength is varied experimentally to assess the trade-off between smoothing and interface accuracy.

% -----------------------------------------------------------
\subsection{Error metrics and evaluation}
Throughout, we evaluate errors against a reference solution $u_{\mathrm{ref}}$ (and
$\mu_{\mathrm{ref}}$ when applicable) computed on a sufficiently fine grid.
For a quantity $v$ (e.g.\ $u$ or $\mu$), the relative $L^2$, relative $L^\infty$, and absolute $L^\infty$ errors are
\[
\mathrm{rel}\,L^2(v)
:= \frac{\|v - v_{\mathrm{ref}}\|_{L^2(Q)}}{\|v_{\mathrm{ref}}\|_{L^2(Q)}}, 
\quad
\mathrm{rel}\, L^\infty(v)
:= \frac{\|v - v_{\mathrm{ref}}\|_{L^\infty(Q)}}{\|v_{\mathrm{ref}}\|_{L^\infty(Q)}}, 
\qquad
L^\infty(v)
:= \|v - v_{\mathrm{ref}}\|_{L^\infty(Q)}.
\]
Unless stated otherwise, the collocation points are chosen as uniform grids for performance evaluation. % (for local MS-RPNN).
%  and the hidden parameters are sampled from $\mathcal U([-R_m,R_m])$ once and fixed.
% \input{results_fun}
\section{Numerical Results}
\subsection{One-dimensional Allen--Cahn equation}\label{sec:AC1d}
We assess the first-order formulation, random-feature design, and accuracy-enhancement strategies on
\begin{equation}\label{eq:ac_1d}
\begin{aligned}
u_t&=\gamma_1u_{xx}-\gamma_2(u^3-u), &&(x,t)\in[-1,1]\times(0,1],\\
u(x,0)&=u_0(x), &&x\in[-1,1],\\
u_x(-1,t)&=u_x(1,t)=0, &&t\in(0,1].
\end{aligned}
\end{equation}
Unless stated otherwise, $\gamma_1=10^{-4}$, $\gamma_2=5$, and $\ell=\sqrt{\gamma_1/\gamma_2}\approx4.47\times10^{-3}$. To assess the robustness of the method across different initial
profiles, we consider $u_0(x)=x^2\cos(\pi x)$,
$u_0(x)=x^2\sin(2\pi x)$, and $u_0(x)=\sin(\pi x)$,
which exhibit different symmetries and spatial structures.

The default configuration uses $N_x=4$, $N_t=1$, and $N_b=10$, with $N_f=2048$ interior, $N_{\mathrm{bc}}=512$ boundary, $N_{\mathrm{ic}}=512$ initial, and $N_{\mathrm{if}}=512$ interface collocation points per local subdomain, as applicable. The initial, boundary, and interface residual multipliers are $100$, $50$, and $50$, respectively. Unless overridden below, each scalar RPNN has $M=400$ features, $\tanh$ activation, separable space--time features, uniform initialization, and $R_\xi=1$.

\paragraph{Impact of activation functions and space--time basis construction.}
We compare seven activation functions, defined in
Table~\ref{tab:activations}, with dense and separable random
features for the sharp-interface AC benchmark
($\ell\approx0.00447$). Table~\ref{tab:basis_act} shows that the smooth activations $\exp$, Gaussian, $\cos$, $\sin$, and $\tanh$ yield similar relative $L^2$ errors, with more noticeable differences in $L^\infty$ error. For these activations, the effect of basis
construction is modest and neither construction is consistently
better. In contrast, the squared ReLU and squared leaky-ReLU
activations are more sensitive to basis construction:
separable features substantially reduce their errors,
although the squared leaky-ReLU retains a large maximum error.
Overall, separable $\tanh$ features achieve the smallest errors
in both metrics and are adopted as the default for subsequent
AC experiments.

\begin{figure} % [htbp]
    \centering
    \includegraphics[width=0.95\linewidth]{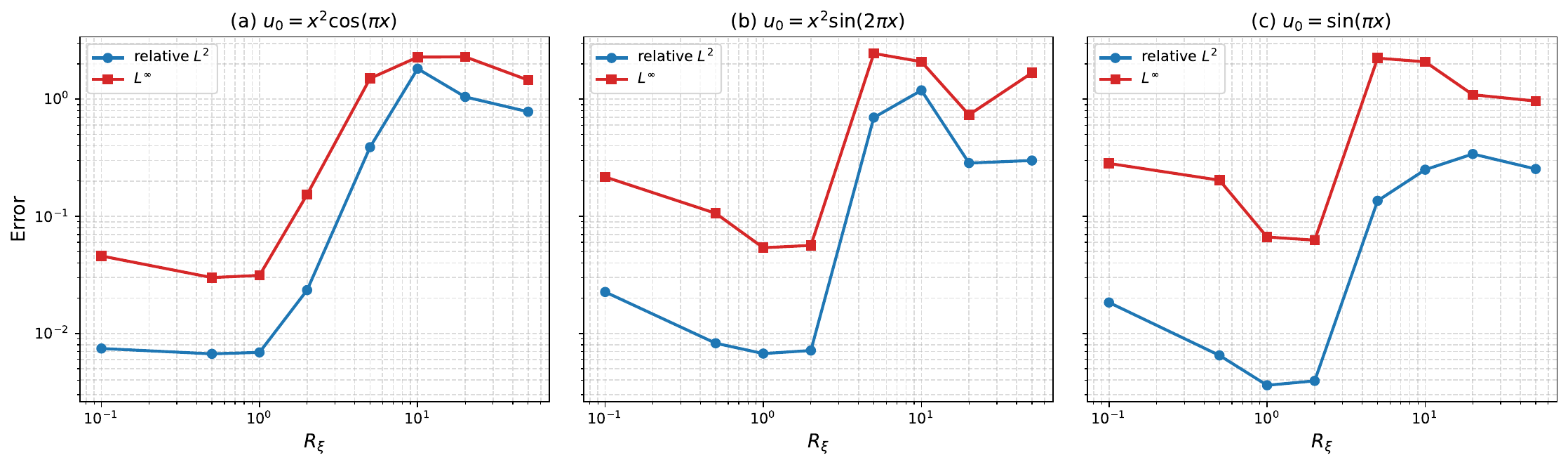}
    \caption{Relative $L^2$ error and $L^\infty$ error with respect to $R_\xi$ for Allen-Cahn equation with three different initial profiles. We set $R_\xi=\frac{a}{\sqrt{3}}$, and $a=[0.1, 0.5, 1, 2, 5, 10, 20, 50]$.}
    \label{fig:ac1d_Rm}
\end{figure}

\paragraph{Impact of random-feature scale $R_\xi$.}
We examine the feature scale using separable tanh features with normal initialization for the sharp-interface AC benchmark with three different initial profiles. As shown in Figure~\ref{fig:ac1d_Rm}, for all three initial profiles, the errors initially decrease with increasing $R_\xi$, then rise sharply, and eventually exhibit a slight drop as $R_\xi$ continues to grow. The best-performing $R_\xi$(yielding the minimum error) varies depending on the specific initial profile, generally falling between 0.5 and 2. Overall, neither too small nor too large values of $R_\xi$ yield the minimum errors. Thus, steeper individual features do not necessarily improve sharp-interface resolution. These results suggest that the feature scale should be selected empirically in relation to the activation, basis construction, and normalized domain, rather than inferred solely from the interface width.

\begin{table}% [htbp]
\centering
\caption{Comparison of dense and separable (sep) bases across activation functions.}
\label{tab:basis_act}
\begin{tabular}{llcc}
\hline
\textbf{Activation} & \textbf{Basis} & \textbf{rel $L^2$} & \textbf{$L^\infty$} \\
\hline
\multirow{2}{*}{exp}
    & dense & $7.45 \times 10^{-3}$ & $5.53 \times 10^{-2}$ \\
    & sep   & $7.42 \times 10^{-3}$ & $5.20 \times 10^{-2}$ \\

\multirow{2}{*}{gauss}
    & dense & $6.89 \times 10^{-3}$ & $6.53 \times 10^{-2}$ \\
    & sep   & $6.92 \times 10^{-3}$ & $7.06 \times 10^{-2}$ \\

\multirow{2}{*}{leakyrelu$^2$}
    & dense & $2.93 \times 10^{-2}$ & $3.92 \times 10^{-1}$ \\
    & sep   & $8.90 \times 10^{-3}$ & $2.48 \times 10^{-1}$ \\

\multirow{2}{*}{relu$^2$}
    & dense & $2.93 \times 10^{-2}$ & $3.92 \times 10^{-1}$ \\
    & sep   & $7.74 \times 10^{-3}$ & $5.98 \times 10^{-2}$ \\

\multirow{2}{*}{cos}
    & dense & $6.93 \times 10^{-3}$ & $6.16 \times 10^{-2}$ \\
    & sep   & $7.00 \times 10^{-3}$ & $6.44 \times 10^{-2}$ \\

\multirow{2}{*}{sin}
    & dense & $6.91 \times 10^{-3}$ & $6.08 \times 10^{-2}$ \\
    & sep   & $6.92 \times 10^{-3}$ & $6.39 \times 10^{-2}$ \\

\multirow{2}{*}{tanh}
    & dense & $6.80 \times 10^{-3}$ & $6.04 \times 10^{-2}$ \\
    & sep   & $\mathbf{6.74 \times 10^{-3}}$ & $\mathbf{3.93 \times 10^{-2}}$\\
\hline
\end{tabular}
\end{table}

% -----------------------------------------------------

% ------------------------------------------------
\begin{table}[htbp]
\centering
\caption{Architectures for the neural networks with $\tanh$ activation.
For FO-RPNN, $u$ and $q$ use independent networks.
$N_p$ denotes the total number of trainable parameters per
space--time subdomain; for RPNNs, only output coefficients are trainable.}
\label{tab:ac1d_structure}
\begin{tabular}{ccccc}
\hline
Model & Architecture & Initialization & Basis type & $N_p$ \\
\hline
PINN-1
& $2 \times 200 \times 1$ & Xavier & -- & 801 \\
PINN-2
& $2 \times 19 \times 19 \times 19 \times 1$ & Xavier & -- & 837 \\
PINN-3
& $2 \times 128 \times 128 \times 128 \times 1$ & Xavier & -- & 33537 \\
\hline
RPNN-d+n & $2 \times 800 \times 1$ & normal & dense & 800 \\
RPNN-d+u & $2 \times 800 \times 1$ & uniform & dense & 800 \\
RPNN-s+n & $2 \times 800 \times 1$ & normal & sep & 800 \\
RPNN-s+u & $2 \times 800 \times 1$ & uniform & sep & 800 \\
\hline
FO-RPNN-d+n & $2 \times 400 \times 1$ for each of $u,q$ & normal & dense & 800 \\
FO-RPNN-d+u & $2 \times 400 \times 1$ for each of $u,q$ & uniform & dense & 800 \\
FO-RPNN-s+n & $2 \times 400 \times 1$ for each of $u,q$ & normal & sep & 800 \\
FO-RPNN-s+u & $2 \times 400 \times 1$ for each of $u,q$ & uniform & sep & 800 \\
\hline
\end{tabular}
\end{table}
 
\paragraph{Impact of initialization and basis construction.}
Table~\ref{tab:ac1d_err_rpnn_vs_pinn} shows that separable
features consistently reduce $L^\infty$ errors relative to
dense features across both RPNN formulations, initialization
distributions, and all three initial profiles. Their advantage
in relative $L^2$ error is less uniform, indicating that basis
construction has a more pronounced effect on maximum error.
Normal initialization generally yields smaller errors for
the sine profiles, whereas neither distribution consistently
dominates for the cosine profile. Separable features with
normal initialization achieve the best results in every
column within the FO-RPNN group and in most columns within
the vanilla RPNN group.

\paragraph{Comparison of PINN, vanilla RPNN, and FO-RPNN.}
We compare PINN and vanilla RPNN applied to the original second-order AC equation~\eqref{eq:ac_1d} with RPNN applied to its equivalent first-order system, denoted by FO-RPNN. Introducing the auxiliary field $q=u_x$ gives
\begin{equation}\label{eq:ac_fo_benchmark}
\begin{aligned}
u_t-\gamma_1q_x+\gamma_2(u^3-u)&=0,\\
q-u_x&=0,
\end{aligned}
\end{equation}
subject to
\[
u(x,0)=u_0(x),
\qquad q(-1,t)=q(1,t)=0.
\]
FO-RPNN approximates $u$ and $q$ using independent random
feature expansions, so evaluating its PDE residuals requires
only first derivatives. The scalar formulation evaluates $u_{xx}$, therefore, the first-order reformulation reduces the derivative order required in residual evaluation, while introducing an additional field and a compatibility residual. The neural network architectures are summarized in Table~\ref{tab:ac1d_structure}. The model performances are summarized in Table~\ref{tab:ac1d_err_rpnn_vs_pinn}. Both RPNN formulations achieve smaller errors than all tested PINN baselines, including PINN-3, despite its substantially larger number of trainable parameters, with reductions of approximately one or two orders of magnitude in many comparisons. The experiments demonstrate the numerical viability of first-order residual regression for this nonlinear AC benchmark. With the same total output-coefficient budget, vanilla RPNN generally outperforms FO-RPNN. Their relative $L^2$ errors remain of comparable magnitude, but FO-RPNN produces larger $L^\infty$ errors, particularly for the sine profiles.This comparison also reflects the allocation of approximation capacity: vanilla RPNN assigns all $800$ coefficients to $u$, while FO-RPNN assigns $400$ each to $u$ and $q$. Thus, the first-order formulation maintains the overall level of accuracy, although reducing the derivative order alone does not consistently improve accuracy under this fixed coefficient budget. 

% ---------------------------------------------
\begin{table}[htbp]
\centering
\small
\caption{Error summary of PINNs, scalar vanilla RPNN and the first-order RPNN (FO-RPNN) for Allen--Cahn equation. The errors are relative $L^2$ and absolute $L^\infty$ errors with respect to the Crank--Nicolson finite difference method.} 
\label{tab:ac1d_err_rpnn_vs_pinn} % \label{tab:first_order_ac_errors}
\resizebox{\linewidth}{!}{%
\begin{tabular}{lcccccc}
\hline
& \multicolumn{2}{c}{$u_0 = x^2\cos(\pi x)$} 
& \multicolumn{2}{c}{$u_0 = x^2\sin(2\pi x)$} 
& \multicolumn{2}{c}{$u_0 = \sin(\pi x)$} \\
\cline{2-7}
\textbf{Model} 
& rel$L^2$ & $L^\infty$ 
& rel$L^2$ & $L^\infty$ 
& rel$L^2$ & $L^\infty$ \\
\hline
PINN-1 & 
$2.71\times10^{-1}$ & $1.07\times10^{0}$ & $2.26\times10^{-1}$ & $9.51\times10^{-1}$ & $1.27\times10^{-1}$ & $7.60\times10^{-1}$ \\
PINN-2 & 
$1.03\times10^{-1}$ & $7.78\times10^{-1}$ & $1.30\times10^{-1}$ & $9.03\times10^{-1}$ & $6.67\times10^{-2}$ & $6.11\times10^{-1}$ \\
PINN-3 & 
$\mathbf{8.85\times10^{-2}}$ & $\mathbf{7.62\times10^{-1}}$ & $\mathbf{1.17\times10^{-1}}$ & $\mathbf{8.42\times10^{-1}}$ & $\mathbf{6.65\times10^{-2}}$ & $\mathbf{5.99\times10^{-1}}$ \\
\hline
RPNN-d+n 
& $6.76\times10^{-3}$ & $5.76\times10^{-2}$ 
& $\mathbf{6.51\times10^{-3}}$ & $1.15\times10^{-1}$ 
& $3.78\times10^{-3}$ & $1.01\times10^{-1}$ \\
RPNN-d+u 
& $6.80\times10^{-3}$ & $6.18\times10^{-2}$ 
& $7.25\times10^{-3}$ & $1.14\times10^{-1}$ 
& $5.62\times10^{-3}$ & $1.85\times10^{-1}$ \\
RPNN-s+n  
& $6.89\times10^{-3}$ & $\mathbf{3.13\times10^{-2}}$ 
& $6.73\times10^{-3}$ & $\mathbf{5.39\times10^{-2}}$ 
& $\mathbf{3.61\times10^{-3}}$ & $\mathbf{6.66\times10^{-2}}$ \\
RPNN-s+u 
& $\mathbf{6.68\times10^{-3}}$ & $3.21\times10^{-2}$ 
& $7.19\times10^{-3}$ & $6.64\times10^{-2}$ 
& $4.54\times10^{-3}$ & $1.21\times10^{-1}$ \\
\midrule
FO-RPNN-d+n & 6.87$\times10^{-3}$ & 6.76$\times10^{-2}$ & 8.23$\times10^{-3}$ & 1.51$\times10^{-1}$ & 6.16$\times10^{-3}$ & 1.94$\times10^{-1}$ \\
FO-RPNN-d+u & 6.78$\times10^{-3}$ & 6.51$\times10^{-2}$ & 7.65$\times10^{-3}$ & 1.37$\times10^{-1}$ & 5.98$\times10^{-3}$ & 2.01$\times10^{-1}$ \\
FO-RPNN-s+n & $\mathbf{6.76\times10^{-3}}$ & $\mathbf{4.27\times10^{-2}}$ & $\mathbf{7.30\times10^{-3}}$ & $\mathbf{8.64\times10^{-2}}$ & $\mathbf{4.99\times10^{-3}}$ & $\mathbf{1.59\times10^{-1}}$ \\
FO-RPNN-s+u & 6.82$\times10^{-3}$ & 5.10$\times10^{-2}$ & 7.54$\times10^{-3}$ & 1.05$\times10^{-1}$ & 5.17$\times10^{-3}$ & 1.66$\times10^{-1}$ \\
\bottomrule
\end{tabular}}
\end{table}

% --------------------------------------------------------
\begin{figure}[htbp]
    \centering
    \includegraphics[width=0.95\linewidth]{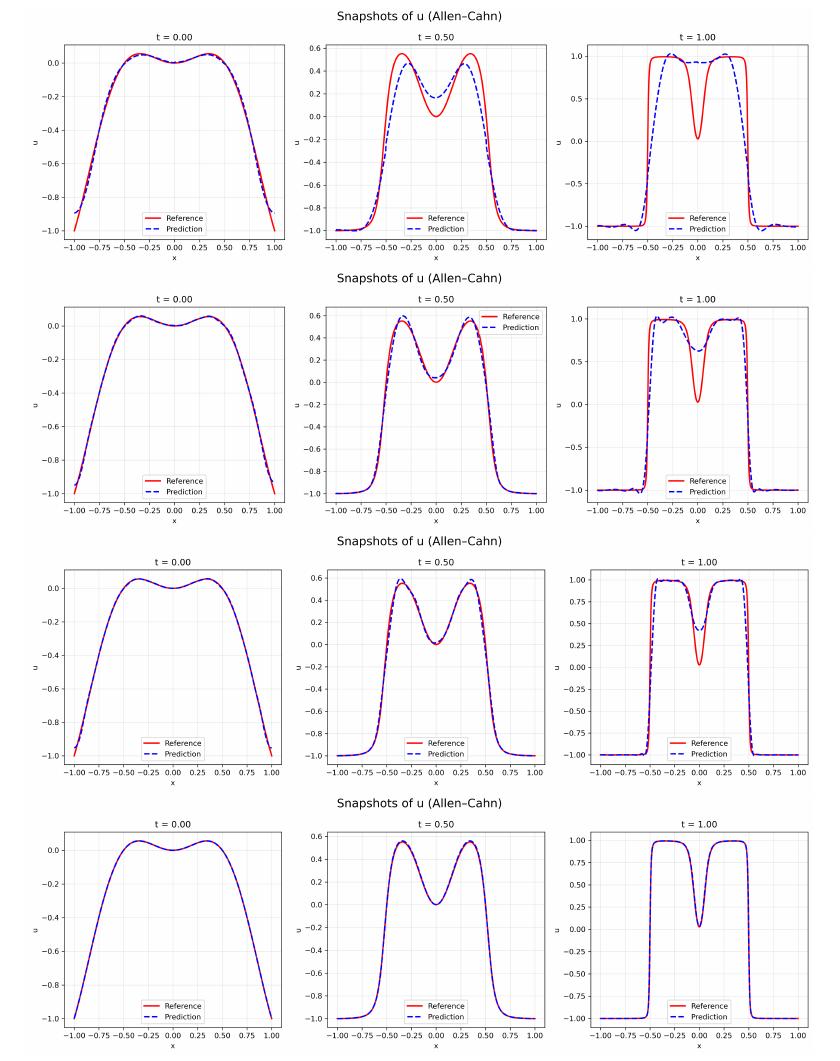}
    \caption{Snapshots at $t=0, 0.5, 1.0$ for the Allen--Cahn equation with $u_0=x^2\cos(\pi x)$: first row: PINN-1; second row: PINN-2; third row: PINN-3; fourth row: scalar vanilla RPNN using separable bases, $\tanh$ activation, and normally distributed random features with $R_\xi=1$ and $M=800$.}
    \label{fig:ac1d_rpnn_vs_pinn}
\end{figure}

\paragraph{Adaptive sampling, reweighting, and regularization.}
Table~\ref{tab:ac1d_tanhsep} compares the baseline with the enhancement strategies defined in Section~\ref{subsec:adap_samp}, using $M=800$, separable $\tanh$ features, and normal initialization. Sampling-based variants enrich the collocation set after an initial solve and restart from its output coefficients; reweighting changes residual contributions without adding points. The listed local/global variants use the corresponding local/global selection pools.

Table~\ref{tab:ac1d_tanhsep} shows that adaptive sampling primarily improves maximum-error control, with larger relative gains for the sine initial profiles. RAR and combined EDRAS+RAR provide the strongest overall improvements: global RAR performs best for $\sin(\pi x)$, while global EDRAS+RAR gives the smallest errors for $x^2\sin(2\pi x)$. However, global selection does not uniformly outperform local selection, and no strategy minimizes every metric.

Grad-ADS improves both error metrics across all three profiles, although its relative $L^2$ gain for the cosine profile is small. EDRAS and residual reweighting consistently reduce maximum errors but do not always improve relative $L^2$ accuracy. Derivative penalties are less reliable: the first-derivative penalty improves the cosine and $\sin(\pi x)$ cases but worsens both errors for
$x^2\sin(2\pi x)$, while the second-derivative penalty degrades accuracy throughout.  

Overall, collocation enrichment is more effective than residual
reweighting and derivative regularization in these experiments,
with benefits most evident in maximum-error reduction;
the best strategy nevertheless depends on the initial profile
and error metric. These findings confirm that adaptive residual- and energy- driven enrichment are promising and cost-effective strategies for enhancing the robustness and approximation power of PDE solvers based on random features.

\begin{table}[htbp]
\centering
\caption{Comparison of different adaptive strategies for the Allen--Cahn equation with three initial conditions using $\tanh$ activation, separable bases, and normal initialization. The suffixes l and g denote local and global selection, respectively.}
\label{tab:ac1d_tanhsep}
\resizebox{\linewidth}{!}{%
\begin{tabular}{lcccccc}
\hline
& \multicolumn{2}{c}{$u_0 = x^2\cos(\pi x)$} 
& \multicolumn{2}{c}{$u_0 = x^2\sin(2\pi x)$} 
& \multicolumn{2}{c}{$u_0 = \sin(\pi x)$} \\
\cline{2-7}
\textbf{Method} 
& rel$L^2$ & $L^\infty$ 
& rel$L^2$ & $L^\infty$ 
& rel$L^2$ & $L^\infty$ \\
\hline
vanilla 
& $6.75\times 10^{-3}$ & $3.93\times 10^{-2}$ 
& $7.50\times 10^{-3}$ & $7.73\times 10^{-2}$ 
& $5.25\times 10^{-3}$ & $1.54\times 10^{-1}$ \\

EDRAS-l
& $6.88\times10^{-3}$ & $3.28\times10^{-2}$ 
& $6.80\times10^{-3}$ & $5.31\times10^{-2}$ 
& $3.37\times10^{-3}$ & $5.15\times10^{-2}$ \\

EDRAS-g
& $6.58\times10^{-3}$ & $2.61\times10^{-2}$ 
& $7.08\times10^{-3}$ & $5.36\times10^{-2}$ 
& $3.50\times10^{-3}$ & $5.37\times10^{-2}$ \\

RAR-l
& $6.65\times10^{-3}$ & $2.10\times10^{-2}$ 
& $6.14\times10^{-3}$ & $3.79\times10^{-2}$ 
& $2.64\times10^{-3}$ & $3.63\times10^{-2}$ \\

RAR-g 
& $6.59\times10^{-3}$ & $2.11\times10^{-2}$ 
& $5.94\times10^{-3}$ & $3.00\times10^{-2}$ 
& $1.90\times10^{-3}$ & $2.17\times10^{-2}$ \\

EDRAS+RAR-l
& $6.64\times10^{-3}$ & $2.18\times10^{-2}$ 
& $6.32\times 10^{-3}$ & $5.13\times 10^{-2}$ 
& $2.62\times10^{-3}$ & $3.65\times10^{-2}$ \\

EDRAS+RAR-g 
& $6.63\times 10^{-3}$ & $2.18\times 10^{-2}$ 
& $5.54\times10^{-3}$ & $2.17\times10^{-2}$ 
& $1.93 \times10^{-3}$ & $ 2.19\times10^{-2}$ \\

Reweighting 
& $6.87\times10^{-3}$ & $3.05\times10^{-2}$ 
& $7.13\times10^{-3}$ & $6.11\times10^{-2}$ 
& $3.82\times10^{-3}$ & $5.99\times10^{-2}$ \\

Grad-ADS
& $6.73\times10^{-3}$ & $2.82\times10^{-2}$ 
& $6.31\times10^{-3}$ & $3.55\times10^{-2}$ 
& $3.02\times10^{-3}$ & $6.46\times10^{-2}$ \\

$u_x$ penalty 
& $6.56\times10^{-3}$ & $2.69\times10^{-2}$ 
& $7.69\times10^{-3}$ & $1.89\times10^{-1}$ 
& $3.37\times10^{-3}$ & $8.89\times10^{-2}$ \\

$u_{xx}$ penalty 
& $1.48\times10^{-2}$ & $2.06\times10^{-1}$ 
& $2.14\times10^{-2}$ & $3.07\times10^{-1}$ 
& $1.56\times10^{-2}$ & $3.56\times10^{-1}$ \\
\hline
\end{tabular}}
\end{table}

% ------------------------------------------------------------------------
\paragraph{Comparison with the literature.} 
We compare our RPNN framework against several state-of-the-art PINN-based methods for solving the Allen–Cahn equation \eqref{eq:AC}, including baseline PINNs~\cite{raissi2019physics}, adaptive PINNs~\cite{wight2020solving}, bc-PINNs~\cite{mattey2022novel}, SA-PINNs~\cite{mcclenny2023self}, PINNs-WE~\cite{liu2024discontinuity}, and Residuals-RAE-PINNs~\cite{zhang2026priori}.

For completeness, we summarize the network architectures employed in these prior works. The baseline PINN, SA-PINN, PINN-WE, and Residuals-RAE-PINN all adopt a fully connected neural network with 2 hidden layers of 128 neurons each, as reported in \cite{zhang2026priori}. Adaptive PINN uses 4 hidden layers with 128 neurons each, while bc-PINN employs 4 hidden layers with 200 neurons each. In contrast, for each subdomain, our RPNN uses a single hidden layer with randomly fixed weights and biases, and determines only 800 output coefficients through least-squares optimization, yielding a substantially smaller model that can be trained more efficiently. We use $N_x=4, N_t=1, N_b=10$. 

Table~\ref{tab:ac_lit_comparison} summarizes the numerical accuracy of these methods alongside our RPNN under various problem settings. Direct comparisons must be interpreted with caution because the reported results correspond to different interface widths \(\ell\) and initial conditions, which significantly affect problem stiffness and achievable accuracy. Nevertheless, our RPNN and FO-RPNN consistently achieve the lowest relative \(L^2\) errors across all tested scenarios. For the moderately sharp case ($\ell\approx1.83\times10^{-2}$), our methods yield errors roughly 2.3 times lower than the best-performing baseline PINN variant (Residuals-RAE-PINN), while outperforming other methods (e.g., SA-PINN, PINN-WE, and baseline PINN) by approximately one to two orders of magnitude. Under the standard sharp-interface setting ($\ell\approx4.47\times10^{-3}$), RPNN and FO-RPNN reduce the error by a factor of 2.5 to 3.5 compared to bc-PINN and adaptive PINN, respectively. These results demonstrate the competitive advantage of the random projection architecture in resolving sharp transition layers with minimal network complexity.

\begin{table}[htbp]
\small
\centering
\caption{Numerical errors for Allen--Cahn equation across different methods and setups.}
\label{tab:ac_lit_comparison}
\resizebox{\linewidth}{!}{%
\begin{tabular}{lcccc}
\hline
Method & Parameters & $\ell$ & Initial condition & rel $L^2$ error \\
\hline
Baseline PINN~\cite{raissi2019physics} & $\gamma_1=10^{-3}, \gamma_2=3$ & $1.83\times10^{-2}$ & $x^2\sin(2\pi x)$ & $1.26\times10^{-1}$ \\
SA-PINN~\cite{mcclenny2023self} & $\gamma_1=10^{-3}, \gamma_2=3$ & $1.83\times10^{-2}$ & $x^2\sin(2\pi x)$ & $6.54\times10^{-2}$ \\
PINN-WE~\cite{liu2024discontinuity} & $\gamma_1=10^{-3}, \gamma_2=3$ & $1.83\times10^{-2}$ & $x^2\sin(2\pi x)$ & $5.60\times10^{-2}$ \\
Residuals-RAE-PINN~\cite{zhang2026priori} & $\gamma_1=10^{-3}, \gamma_2=3$ & $1.83\times10^{-2}$ & $x^2\sin(2\pi x)$ & $8.09\times10^{-3}$ \\
\textbf{RPNN (ours)} & $\gamma_1=10^{-3}, \gamma_2=3$ & $1.83\times10^{-2}$ & $x^2\sin(2\pi x)$ & $\mathbf{3.46\times10^{-3}}$ \\
\textbf{FO-RPNN (ours)} & $\gamma_1=10^{-3}, \gamma_2=3$ & $1.83\times10^{-2}$ & $x^2\sin(2\pi x)$ & $\mathbf{3.65\times10^{-3}}$ \\
\hline
Adaptive PINN~\cite{wight2020solving} & $\gamma_1=10^{-4}, \gamma_2=5$ & $4.47\times10^{-3}$ & $x^2\cos(\pi x)$ & $2.33\times10^{-2}$ \\
bc-PINN~\cite{mattey2022novel} & $\gamma_1=10^{-4}, \gamma_2=5$ & $4.47\times10^{-3}$ & $x^2\cos(\pi x)$ & $1.68\times10^{-2}$ \\
\textbf{RPNN (ours)} & $\gamma_1=10^{-4}, \gamma_2=5$ & $4.47\times10^{-3}$ & $x^2\cos(\pi x)$ & $\mathbf{6.68\times10^{-3}}$ \\
\textbf{FO-RPNN (ours)} & $\gamma_1=10^{-4}, \gamma_2=5$ & $4.47\times10^{-3}$ & $x^2\cos(\pi x)$ &  $\mathbf{6.76\times 10^{-3}}$\\
\hline
\end{tabular}}
\end{table}
\subsection{One-dimensional Cahn--Hilliard equation}\label{sec:results_ch}
We consider the unit-mobility Cahn--Hilliard equation with homogeneous Neumann conditions on both the order parameter and the physical chemical potential:
\begin{equation}\label{eq:CH_1d}
\begin{aligned}
u_t-\mu_{xx}&=0, &\mu&=\gamma_2(u^3-u)-\gamma_1u_{xx},\\
u(x,0)&=u_0(x), &u_x&=\mu_x=0\quad\text{at the spatial endpoints}.
\end{aligned}
\end{equation}
The parameter $\ell := \sqrt{\frac{\gamma_1}{\gamma_2}}$ is an indicator of the interface width: smaller $\ell$ corresponds to sharper interfaces that are more challenging to resolve accurately. The primary benchmark (Case A) uses $x\in[-1,1]$, $t\in[0,1]$, $u_0(x)=-\cos(2\pi x)$, $\gamma_1=10^{-6}$, and $\gamma_2=10^{-2}$, giving $\ell=\sqrt{\gamma_1/\gamma_2}=0.01$. Case B uses $x\in[0,2\pi]$, $t\in[0,4]$, $u_0(x)=0.5\cos x$, $\gamma_1=0.03$, and $\gamma_2=1$, giving $\ell=\sqrt{0.03}$. These cases test distinct initial profiles and interface scales.

Unless stated otherwise, the mixed two-field experiments use $N_x=8$, $N_t=1$, and $N_b=40$, with $N_f=1024$ and $N_{\mathrm{bc}}=N_{\mathrm{ic}}=N_{\mathrm{if}}=256$ per local subdomain, as applicable. Initial, boundary, and interface residual multipliers are $100$, $50$, and $50$. Independent single-scale $\tanh$ features approximate the two fields with $M_u=800$, $M_{\mu_s}=120$, and feature scales equal to one. The evolution and chemical-potential residual factors are $1$ and $100$, unless varied explicitly below.

\begin{figure}[htbp]
    \centering
    \includegraphics[width=0.95\linewidth]{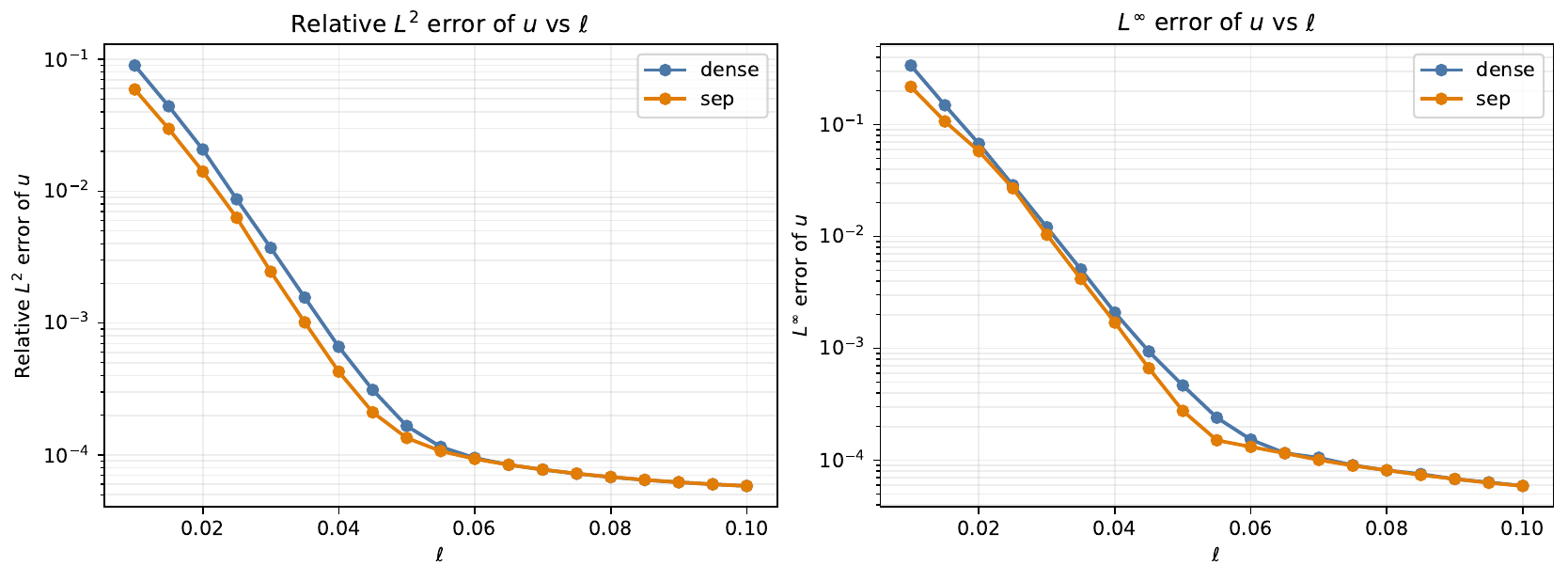}
    \caption{Errors in $u$ versus interface width $\ell$ for the fixed shared-basis configuration: relative $L^2$ (left) and absolute $L^\infty$ (right). }
    \label{fig:relL2_u_l}
\end{figure}

\paragraph{Effect of interface width $\ell$.}
We investigate the effect of interface sharpness by varying $\ell$, using $\tanh$ activation and a shared random basis of $M=400$ features with $R_\xi=1$ for $u$ and $\mu$ in \eqref{eq:CH_1d}. Figure~\ref{fig:relL2_u_l} reports the relative $L^2$ error and $L^{\infty}$ error in $u$ for each $\ell$ for dense and separable bases. The results show that accuracy deteriorates as $\ell$ decreases, as expected, since thinner interfaces produce sharper gradients and are more difficult to resolve. In particular, the error increases rapidly once $\ell\lesssim 0.05$. 

\paragraph{Chemical-potential rescaling.}
A key challenge in solving the Cahn–Hilliard equation with neural network solvers is the computation of fourth-order spatial derivatives, which introduces substantial memory and computational costs due to automatic differentiation. To reduce the order of spatial derivatives, we introduce an auxiliary variable $\mu$ and consider three equivalent scalings, all with the same boundary and initial conditions. To separate the physical field from its numerical representation, write $\mu=s\mu_s$, $s>0$. The mixed equations become
\begin{equation}\label{eq:ch_mixed_scaled_results}
u_t-s(\mu_s)_{xx}=0,
\qquad
\mu_s-\frac{\gamma_2}{s}(u^3-u)+\frac{\gamma_1}{s}u_{xx}=0.
\end{equation}
The three cases correspond to $s=1$ (Case 1), $s=\gamma_1$ (Case 2), and $s=\gamma_2$ (Case 3). We call this mixed two-field solver vanilla RPNN to distinguish it from the four-field first-order solver; it already avoids direct fourth derivatives but still evaluates second derivatives. The three scalings preserve the exact physical solution, but change the auxiliary-field magnitude.

\begin{figure} %[htbp]
    \centering
    \includegraphics[width=0.9\linewidth]{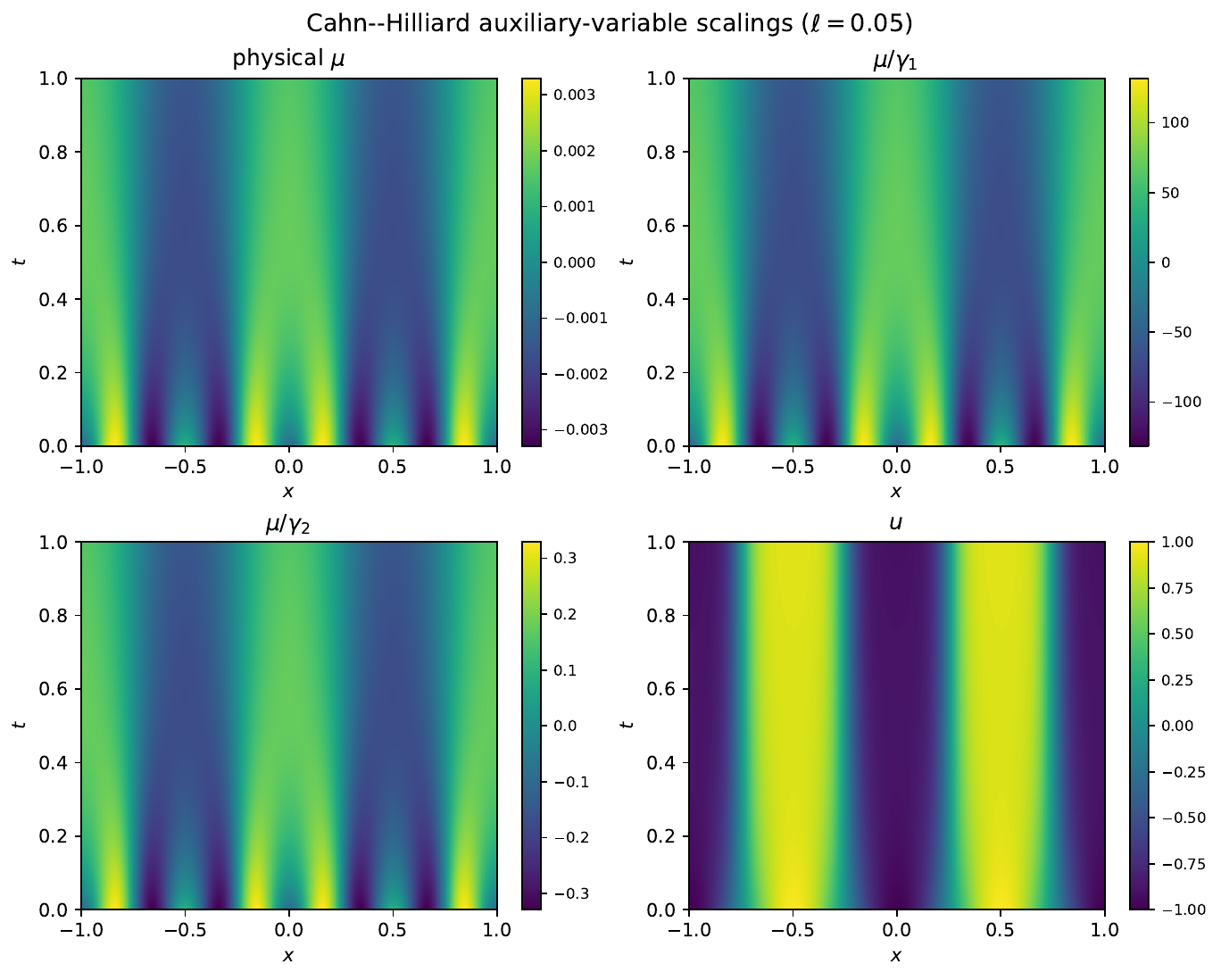}
    \caption{Auxiliary variable $\mu_s$ under the three scalings. Case 1 (top left) and Case 2 (top right) both show severe imbalance between $u$ and $\mu_s$; Case 3 (bottom left) shows a well-balanced relation between $u$ and $\mu_s$; the bottom right panel shows $u$.} 
    \label{fig:3mu}
\end{figure}

Although the three formulations are mathematically equivalent, they produce markedly different magnitude balances between \(u\) and the auxiliary variable \(\mu_s\), as shown in Figure~\ref{fig:3mu}. Cases 1 and 2 exhibit severe scale imbalances: \(\mu_s\) is much smaller than \(u\) in Case 1 and much larger in Case 2. In contrast, Case 3 gives \(\mu_s=(u^3-u)-\ell^2u_{xx}\), bringing the two fields to comparable magnitudes. This scaling normalizes the nonlinear coefficient and provides a more balanced representation for their simultaneous approximation. The resulting balance or imbalance between the field magnitudes can substantially affect the conditioning of the least-squares problem when \(u\) and \(\mu_s\) are approximated simultaneously.

Table~\ref{tab:ch_mu_scaling} summarizes the errors for $\gamma_1=2.5\times 10^{-5}, \ell=0.05$ under the three equivalent scalings cases. In this comparison, we employ independent separable random bases for \(u\) and \(\mu\), so as to circumvent potential performance degradation that may arise from the limited expressivity or structural restrictions of a single shared neural network. The results reveal that the choice of scaling has a dramatic effect on both accuracy and stability of the \((u,\mu)\) approximation. For Case~1 and Case~2, the relative \(L^2\) errors for \(u\) are nearly \(1\), indicating a complete failure of the approximation. In contrast, Case~3 yields relative \(L^2\) errors on the order of \(10^{-4}\) for both \(u\) and \(\mu\), which is orders of magnitude more accurate. Notably, within the successful Case 3, the dense basis yields slightly smaller errors than the separable basis. 

This dramatic difference in accuracy is directly tied to the system conditioning. In Case~1 and Case~2, the magnitudes of the two fields differ by several orders, leading to an ill-conditioned least-squares system. In Case~3, however, the coefficients ensure \(u\) and \(\mu\) are of comparable magnitude, yielding a better-balanced residual. This is explicitly supported by the global-Jacobian conditioning analysis in Table~\ref{tab:ch_global_jacobian_conditioning}: the effective condition number \(\kappa_{\text{eff}}\) for Case 3 is lower than those of Cases 1 and 2. Moreover, Table~\ref{tab:ch_global_jacobian_conditioning} shows that the dense basis consistently provides a lower \(\kappa_{\text{eff}}\) than the separable basis across all cases (with a Dense/Sep ratio of 0.385 in Case 3), which explains the marginal accuracy advantage of the dense representation observed in Table~\ref{tab:ch_mu_scaling}. These observations indicate that the \(\gamma_2\)-scaling (Case~3) is the preferred formulation for the auxiliary-variable approach in the Cahn--Hilliard equation, as it ensures both high accuracy and numerical robustness. In the following experiments, we always choose Case~3 to define the auxiliary variable.

\paragraph{Independent bases and residual weighting.} When the proper scaled auxiliary variable $\mu$ is introduced, a natural question arises: \textit{should we use a shared random basis for both $u$ and $\mu$, or independent bases?} Independent bases allow feature widths and scales to be tailored to the approximation  requirements of each field. The order parameter $u$ and the chemical potential $\mu$ often differ in smoothness and magnitude, particularly in sharp-interface regimes. Using two distinct sets of random features therefore provides the flexibility to choose different values for $R_{\xi_u}$ and $R_{\xi_{\mu}}$ independently.  To investigate whether $u$ and $\mu_s$ benefit from shared or independent random bases, we conducted a series of numerical experiments, with results summarized in Table~\ref{tab:ch_shared_basis} which distinguishes shared features from independent feature realizations; here, independence refers to separate feature realizations for $u$ and $\mu$, not space--time separability. At equal widths $M_u=M_{\mu_s}=200$, independent bases reduce both error metrics for both fields. The independent $800+120$ configuration also outperforms the shared $800+800$ configuration. These results demonstrate the potential of field-specific feature allocation: the chemical potential can be represented with substantially fewer features while improving the accuracy of both fields in this benchmark. Unless otherwise stated, subsequent experiments therefore use independent random bases for $u$ and $\mu_s$.
   
Table~\ref{tab:ch1d_err_rpnn_ablation} shows that residual scaling remains important after variable rescaling. Increasing the reported chemical-potential residual factor from $1$ to $100$ improves every listed RPNN configuration. Variable scaling and equation weighting are distinct design choices.

% ---------------------------------------------------------
% ----------------------------------------------
\begin{table} %[htbp]
\centering
\caption{Scaling comparison at $\ell=0.05$ for Cahn--Hilliard equation.}
\label{tab:ch_mu_scaling}
\begin{tabular}{clcccc}
\toprule
Case & Basis & \multicolumn{2}{c}{$u$} & \multicolumn{2}{c}{$\mu_s$} \\
\cmidrule(lr){3-4} \cmidrule(lr){5-6}
 & & rel $L^2$ & $L^\infty$ & rel $L^2$ &  $L^\infty$ \\
\midrule
case 1  & dense  & $9.998\times10^{-1}$ & $9.934\times10^{-1}$ & $1.000\times10^{0}$ &  $3.293\times10^{-3}$ \\
case 1 & sep    & $9.999\times10^{-1}$ & $9.942\times10^{-1}$ & $1.000\times10^{0}$ &  $3.293\times10^{-3}$ \\
case 2 ($\gamma_1$) & dense & $1.000\times10^{0}$ &  $1.000\times10^{0}$ & $1.000\times10^{0}$ &  $1.317\times10^{2}$ \\
case 2 ($\gamma_1$) & sep   & $1.000\times10^{0}$ &  $1.000\times10^{0}$ & $9.999\times10^{-1}$ &  $1.317\times10^{2}$ \\
case 3 ($\gamma_2$) & dense & $7.435\times10^{-5}$ & $1.980\times10^{-4}$ & $3.103\times10^{-4}$ &  $1.906\times10^{-4}$ \\
case 3 ($\gamma_2$) & sep   & $8.373\times10^{-5}$ & $2.015\times10^{-4}$ & $4.163\times10^{-4}$ & $3.081\times10^{-4}$ \\
\bottomrule
\end{tabular}
\end{table}

% ------------------------------------------
\begin{table} %[htbp]
\centering
\caption{Global-Jacobian conditioning in the first sequential time slab under the implemented loss weights $w_{u,\mathrm{PDE}}=1$, $w_{\mu,\mathrm{PDE}}=100$, $w_{\mathrm{IC}}=100$, and $w_{\mathrm{BC}}=w_{\mathrm{face}}=50$ for $\ell=0.05$. The effective condition number $\kappa_{\mathrm{eff}}$ excludes singular values below the numerical-rank threshold. All cases are numerically rank deficient; consequently their full-column-rank condition numbers are infinite. $r_n$ means numerical rank.}
\label{tab:ch_global_jacobian_conditioning}
\begin{tabular}{lrrrrrr}
\toprule
$\mu$ formulation
& Dense $\kappa_{\mathrm{eff}}$
& Sep $\kappa_{\mathrm{eff}}$
& Dense/Sep
&Dense $r_n$ 
& Sep $r_n$ 
& Total columns \\
\midrule
Case 1
& $4.074\times10^{10}$ & $5.077\times10^{10}$ & $0.802$
& 5192 & 4834 & 7360 \\
Case 2
& $1.451\times10^{11}$ & $1.460\times10^{11}$ & $0.994$
& 5110 & 4725 & 7360 \\
Case 3
& $\mathbf{1.950\times10^{10}}$ & $5.059\times10^{10}$
& $\mathbf{0.385}$ & 5281 & 4888 & 7360 \\
\bottomrule
\end{tabular}
\end{table}

% -----------------------------------------------------------------
\begin{table}% [htbp]
\centering
\caption{Relative and absolute errors for $u$ and $\mu_s$ using shared and independent bases. Dense bases with $\tanh$ activation and uniform initialization are used to generate the random feature functions.} 
\label{tab:ch_shared_basis}
\resizebox{\linewidth}{!}{%
\begin{tabular}{cclcccc}
\toprule
Case & $M_u$ & $M_{\mu_s}$ & \multicolumn{2}{c}{$u$} & \multicolumn{2}{c}{$\mu_s$} \\
\cmidrule(lr){4-5} \cmidrule(lr){6-7}
& & & rel $L^2$ & $L^\infty$ & rel $L^2$ & $L^\infty$ \\
\midrule
shared basis & 800 & 800
& $1.229\times10^{-3}$ & $9.883\times10^{-3}$
& $1.601\times10^{-3}$ & $3.481\times10^{-3}$ \\

independent basis & 800 & 120 &
 $8.638\times10^{-4}$ & $5.299\times10^{-3}$ & $1.024\times10^{-3}$ & $1.656\times10^{-3}$ \\

shared basis & 200 & 200
& $3.829\times10^{-3}$ & $4.006\times10^{-2}$
& $1.265\times10^{-2}$ & $3.213\times10^{-2}$ \\

independent basis & 200 & 200
& $1.327\times10^{-3}$ & $1.23\times10^{-2}$
& $2.332\times10^{-3}$ & $3.703\times10^{-3}$ \\
\bottomrule
\end{tabular}}
\end{table}

\begin{figure}[htbp]
    \centering
    \includegraphics[width=0.99\linewidth]{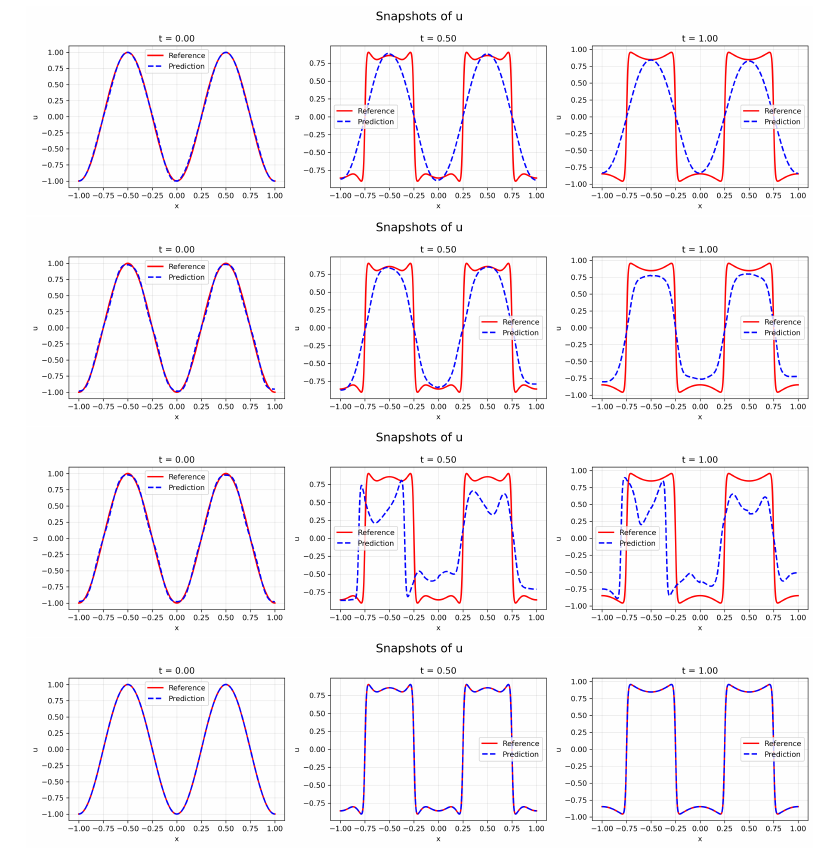}
    \caption{Snapshots at $t=0, 0.5, 1.0$ for the Cahn--Hilliard equation with $u_0=-\cos(2\pi x)$: first row: PINN-1; second row: PINN-2; third row: PINN-3; fourth row: mixed RPNN using dense bases, $\tanh$ activation, and uniformly distributed random features with $R_\xi=1$ and $M=800$.}
    \label{fig:ch_rpnn_pinn_cos_uniform}
\end{figure}

\paragraph{Mixed PINN, mixed RPNN, and first-order RPNN.}
We compare mixed PINNs, mixed RPNNs, and first-order RPNNs (FO-RPNNs) for the Cahn--Hilliard equation. Table~\ref{tab:ch1d_structure} summarizes their architectures and numbers of trainable parameters per space--time subdomain. PINN-1, PINN-2, mixed RPNN, and FO-RPNN have comparable parameter counts (ranging from $800$ to $978$), whereas PINN-3 provides a larger-capacity baseline with $67{,}074$ trainable parameters. All models use the \(\tanh\) activation function. The mixed RPNNs employ independent random feature bases for $u$ and the scaled chemical potential $\mu_s$ to accommodate their different characteristic scales and regularity properties. Solution snapshots at $t=0$, $0.5$, and $1$ are shown in Figures~\ref{fig:ch_rpnn_pinn_cos_uniform} and~\ref{fig:ch_rpnn_pinn_cos_uniform_mu}.

For Case A, Table~\ref{tab:ch1d_err_rpnn_vs_pinn} shows that the PINNs yield high relative $L^2$ errors on the order of \(10^{-1}\) for both fields. Notably, increasing the network size does not improve accuracy; in fact, PINN-3 produces the largest errors among the three PINN architectures, indicating severe optimization difficulties within the same number of training epochs. By contrast, mixed RPNNs with uniformly distributed features dramatically reduce the errors to the \(10^{-3}\)--\(10^{-4}\) range, achieving an improvement of roughly two to three orders of magnitude over the PINNs. Dense and separable uniform features perform comparably well, both yielding errors well below \(10^{-3}\). Among all FO-RPNN results, separable uniform features perform the best, achieving an accuracy of the same order (\(10^{-3}\)) as the best mixed RPNN results.

For Case B, Table~\ref{tab:ch1d_err_rpnn_halfcos} reveals a clearer advantage for the first-order formulation. All mixed RPNN configurations yield nearly identical errors of order \(10^{-4}\), which are already over two orders of magnitude better than the best PINN baselines. FO-RPNNs further reduce these errors to the order of \(10^{-5}\), corresponding to improvements by factors of about $6$ to $7$ for $u$ and over $20$ for $\mu_s$ compared to mixed RPNNs. Their absolute $L^\infty$ errors are similarly reduced by more than an order of magnitude. Furthermore, differences between dense and separable features, as well as between uniform and normal initialization, are negligible within each RPNN formulation for this case.

Overall, the RPNN configurations substantially outperform the PINN baselines in both cases despite using far fewer trainable parameters than PINN-3 and having parameter counts comparable to those of PINN-1 and PINN-2. FO-RPNN provides a further, distinct accuracy improvement in Case B, while Case A highlights the sensitivity of the mixed formulation to feature construction and residual scaling. These findings confirm that optimizing only the output coefficients of fixed random features is highly effective for these problems.

% -----------------------------------------------
\paragraph{Published PINN comparisons.}
In this part, we compare our method with the results in literature. We collect the results from baseline \cite{raissi2019physics}, adaptive PINNs\cite{wight2020solving}, bc-PINN\cite{mattey2022novel}, SA-PINNs\cite{mcclenny2023self}, PINNs-WE\cite{liu2024discontinuity}, and Residuals-RAE-PINNs\cite{zhang2026priori} that using neural networks to solve the Cahn--Hilliard equation given below:
\begin{equation}
\begin{aligned}
u_t &=  \Delta \mu,  \quad  \bx \in [-1, 1], t\in (0, T],\\
\mu &= -\gamma_1 \Delta u + \gamma_2(u^3-u), \quad \bx \in [-1, 1], t\in (0, T],
\end{aligned}
\end{equation}
where periodic boundary condition is utilized and the initial condition as well as equation parameters are listed in Table~\ref{tab:ch_lit_comparison}.

Note that the terminal time is $T=0.5$ for baseline\cite{raissi2019physics}, SA-PINNs\cite{mcclenny2023self}, PINNs-WE\cite{liu2024discontinuity}, Residuals-RAE-PINNs\cite{zhang2026priori}. Meanwhile, the terminal time is $T=1.0$ for adaptive-PINN \cite{wight2020solving} and bc-PINN \cite{mattey2022novel}. As demonstrated in Table~\ref{tab:ch_lit_comparison}, our proposed Mixed RPNN and FO-RPNN consistently achieve the lowest relative \(L^2\) errors on the order of \(10^{-3}\) to \(10^{-4}\), significantly outperforming all the compared PINN baselines and variants. Specifically, under the standard parameter settings (\(\gamma_1=0.02, \gamma_2=1, \ell=0.14\)), our method reduces the error by nearly two orders of magnitude compared to the baseline, SA-PINNs, and PINNs-WE, and outperforms the best-performing PINN variant (Residuals-RAE-PINNs) by over an order of magnitude. It also achieves a 5-fold improvement over bc-PINN. Furthermore, even in the extremely sharp interface case (\(\gamma_1=10^{-6}, \ell=0.01\)), our Mixed RPNN yields an error more than 10 times lower than the adaptive-PINN. These results robustly demonstrate the competitive advantage of the random projection architecture in resolving sharp transition layers across diverse problem setups with minimal network complexity.

% ---------------------------------------------------  
\begin{table}[htbp]
\centering
\caption{Architectures for the neural networks with $\tanh$ activation. Mixed PINNs and Mixed PINNs use two independent networks and FO-PINNs use four independent networks. RPNNs train only output coefficients.}
\label{tab:ch1d_structure}
\resizebox{\linewidth}{!}{%
\begin{tabular}{llll}
\toprule
Model & Architecture or feature widths & Initialization/basis & Total\\
\midrule
PINN-1 & $u$: $2\times200\times1$; $\mu_s$: $2\times30\times1$ & Xavier / dense & 922\\
PINN-2 & $u$: $2\times19\times19\times19\times1$; $\mu_s$: $2\times7\times7\times7\times1$ & Xavier / dense & 978\\
PINN-3 & Each field: $2\times128\times128\times128\times1$ & Xavier / dense & 67074\\
Mixed RPNN & $(M_u,M_{\mu_s})=(800,120)$ & Uniform or normal / dense or separable & 920\\
FO-RPNN & $(M_u,M_q,M_\mu,M_j)=(200,200,200,200)$ & Uniform or normal / dense or separable & 800\\
\bottomrule
\end{tabular}}
\end{table}

% ---------------- case A: Best model comparison (w=100) -----------------------------
\begin{table}[htbp]
\centering
\caption{Error comparison for PINNs, Mixed RPNNs and FO-RPNNs (Case A: \(u(x,0)=-\cos(2\pi x)\)). $w_{\mu_s,\mathrm{PDE}}$ denotes the reported chemical-potential residual scaling factor.}
\label{tab:ch1d_err_rpnn_vs_pinn}
\resizebox{\linewidth}{!}{%
\begin{tabular}{lccccc}
\toprule
Model Structure & \(w_{\mu_s, \text{PDE}}\) & \(\text{rel} L^2_u\) & \(L^{\infty}_u\) & \(\text{rel} L^2_{\mu_s}\) & \(L^{\infty}_{\mu_s}\) \\
\midrule
PINN-1 & 100 & \(4.014\times 10^{-1}\) & \(8.093\times10^{-1}\) & \(3.673\times10^{-1}\) & \(2.281\times10^{-1}\) \\
PINN-2 & 100 & \(3.523\times10^{-1}\) & \(8.006\times10^{-1}\) & \(4.105\times10^{-1}\) & \(2.843\times10^{-1}\) \\
PINN-3 & 100 & \(6.020\times10^{-1}\) & \(1.924\times10^{0}\) & \(7.914\times10^{-1}\) & \(8.949\times10^{-1}\) \\
\midrule
Mixed RPNN-d+u & 100 & \(\mathbf{8.638\times10^{-4}}\) & \(\mathbf{5.299\times10^{-3}}\) & \(\mathbf{1.024\times10^{-3}}\) & \(\mathbf{1.656\times10^{-3}}\) \\
Mixed RPNN-d+n & 100 & \(4.433\times10^{-3}\) & \(4.237\times10^{-2}\) & \(8.925\times10^{-3}\) & \(1.915\times10^{-2}\) \\
Mixed RPNN-s+u & 100 & \(\mathbf{8.896\times10^{-4}}\) & \(\mathbf{4.824\times10^{-3}}\) & \(\mathbf{8.633\times10^{-4}}\) & \(\mathbf{1.389\times10^{-3}}\) \\
Mixed RPNN-s+n & 100 & \(7.249\times 10^{-3}\) & \(5.344\times 10^{-2}\) & \(5.873\times 10^{-3}\) & \(1.002\times 10^{-2}\) \\
\midrule
FO-RPNN-d+n & 100 & \(1.84\times10^{-2}\) & \(1.13\times10^{-1}\) & \(5.32\times10^{-2}\) & \(2.11\times10^{-3}\) \\
FO-RPNN-d+u & 100 & \(9.04\times10^{-3}\) & \(7.57\times10^{-2}\) & \(3.08\times10^{-2}\) & \(8.97\times10^{-4}\) \\
FO-RPNN-s+n & 100 & \(1.90\times10^{-3}\) & \(2.89\times10^{-2}\) & \(1.10\times10^{-2}\) & \(6.83\times10^{-4}\) \\
FO-RPNN-s+u & 100 & \(\mathbf{1.15\times10^{-3}}\) & \(\mathbf{8.10\times10^{-3}}\) & \(\mathbf{4.37\times10^{-3}}\) & \(\mathbf{1.46\times10^{-4}}\) \\
\bottomrule
\end{tabular}}
\end{table}

% ----------------------case B ------------------------
\begin{table}[htbp]
\centering
\caption{Error comparison for different neural network structures (Case B: $u(x,0)=0.5\cos(x)$). $w_{\mu_s,\mathrm{PDE}}$ denotes the reported chemical-potential residual scaling factor.}
\label{tab:ch1d_err_rpnn_halfcos}
\resizebox{\linewidth}{!}{%
\begin{tabular}{lccccc}
\hline
Model Structure 
& $w_{\mu_s,\mathrm{PDE}}$ 
& $\mathrm{rel}\,L^2_u$ 
& $L^\infty_u$ 
& $\mathrm{rel}\,L^2_{\mu_s}$ 
& $L^\infty_{\mu_s}$ \\
\hline
PINN-1
& 100
& $5.907\times10^{-2}$
& $2.785\times10^{-1}$
& $\mathbf{1.714\times10^{-1}}$
& $\mathbf{9.422\times10^{-2}}$ \\

PINN-2
& 100
& $\mathbf{5.057\times10^{-2}}$
& $\mathbf{2.077\times10^{-1}}$
& $2.198\times10^{-1}$
& $9.504\times10^{-2}$ \\

PINN-3
& 100
& $1.428\times10^{0}$
& $1.988\times10^{0}$
& $1.339\times10^{0}$
& $8.814\times10^{-1}$ \\
\hline

Mixed RPNN-d+u
& 100
& $3.561\times10^{-4}$
& $1.706\times10^{-3}$
& $6.388\times10^{-4}$
& $\mathbf{5.726\times10^{-4}}$ \\

Mixed RPNN-d+n
& 100
& $3.560\times10^{-4}$
& $1.706\times10^{-3}$
& $\mathbf{6.381\times10^{-4}}$
& $5.742\times10^{-4}$ \\

Mixed RPNN-s+u
& 100
& $3.561\times10^{-4}$
& $\mathbf{1.706\times10^{-3}}$
& $6.388\times10^{-4}$
& $5.726\times10^{-4}$ \\

Mixed RPNN-s+n
& 100
& $\mathbf{3.560\times10^{-4}}$
& $1.706\times10^{-3}$
& $6.383\times10^{-4}$
& $5.743\times10^{-4}$ \\
\hline
FO-RPNN-d+n & 100 & 5.4155$\times10^{-5}$ & 1.4442$\times10^{-4}$ & 3.1075$\times10^{-5}$ & 2.7411$\times10^{-5}$ \\
FO-RPNN-d+u & 100 & $\mathbf{5.4098\times10^{-5}}$ & 1.4343$\times10^{-4}$ & 3.0979$\times10^{-5}$ & 2.7114$\times10^{-5}$ \\
FO-RPNN-s+n & 100 & 5.4103$\times10^{-5}$ & 1.4211$\times10^{-4}$ & 3.0931$\times10^{-5}$ & 2.7104$\times10^{-5}$ \\
FO-RPNN-s+u & 100 & 5.4102$\times10^{-5}$ & $\mathbf{1.4209\times10^{-4}}$ & $\mathbf{3.0930\times10^{-5}}$ & $\mathbf{2.7103\times10^{-5}}$ \\
\bottomrule
\end{tabular}}
\end{table}

% ---------------- literature comparision --------------------------
\begin{table}[htbp]
\small
\centering
\caption{Numerical errors for Cahn--Hilliard equation across different methods and setups.}
\label{tab:ch_lit_comparison}
\resizebox{\linewidth}{!}{%
\begin{tabular}{lcccc}
\hline
Method & Parameters & $\ell$ & I.C. & rel $L^2$ Error \\
\hline
baseline\cite{raissi2019physics} & $\gamma_1=0.02, \gamma_2=1$ & $0.14$ & $\cos(\pi x)-\exp(-4(\pi x)^2)$ & $3.73\times 10^{-1}$\\
 
SA-PINNs~\cite{mcclenny2023self} & $\gamma_1=0.02, \gamma_2=1$ & $0.14$ & $\cos(\pi x)-\exp(-4(\pi x)^2)$ & $2.90\times 10^{-1}$\\
 
PINNs-WE~\cite{liu2024discontinuity} & $\gamma_1=0.02, \gamma_2=1$ & $0.14$ & $\cos(\pi x)-\exp(-4(\pi x)^2)$  & $4.72\times 10^{-1}$\\
Residuals-RAE-PINNs~\cite{zhang2026priori} & $\gamma_1=0.02, \gamma_2=1$ & $0.14$ & $\cos(\pi x)-\exp(-4(\pi x)^2)$ & $4.83\times 10^{-2}$\\
\textbf{Mixed RPNN (ours)} & $\gamma_1=0.02, \gamma_2=1$ & $0.14$ & $\cos(\pi x)-\exp(-4(\pi x)^2)$ & $\mathbf{3.61\times 10^{-3}}$ \\

\hline
adaptive-PINN \cite{wight2020solving} & $\gamma_1=10^{-6}, \gamma_2=10^{-2}$ & $0.01$ & $-\cos(2\pi x)$ & $9.51 \times 10^{-3}$\\
\textbf{Mixed RPNN (ours)} & $\gamma_1=10^{-6}, \gamma_2=10^{-2}$ & $0.01$ & $-\cos(2\pi x)$ & $\mathbf{8.64\times 10^{-4}}$\\
\hline

bc-PINN \cite{mattey2022novel} & $\gamma_1=0.02, \gamma_2=1$ & $0.14$ &$\cos(\pi x)-\exp(-4(\pi x)^2)$ &$ 1.86\times 10^{-2}$\\
\textbf{Mixed RPNN (ours)} & $\gamma_1=0.02, \gamma_2=1$ & $0.14$ &$\cos(\pi x)-\exp(-4(\pi x)^2)$ & $ \mathbf{3.64\times 10^{-3}}$ \\
\hline
\end{tabular}}
\end{table}

% ---------------------------------------------------

\section{Conclusion}\label{sec:conclusion}

We developed an RPNN least-squares collocation framework for the one-dimensional Allen--Cahn (AC) and Cahn--Hilliard (CH) equations, investigating first-order auxiliary-variable formulations alongside their higher-order counterparts. Fixed hidden-layer parameters reduce training to the determination of output coefficients using nonlinear least-squares solvers, while independent random-feature expansions provide flexibility in approximating the coupled fields. The experiments demonstrate the numerical viability of first-order residual regression for nonlinear AC and CH benchmarks. 

Our experiments identify chemical-potential scaling and field-specific basis construction as important design choices for the CH equation. In the mixed formulation, scaling by the potential coefficient $\gamma_2$ substantially improves accuracy in the reported scaling test, while independent bases permit effective approximation with different feature counts for the two fields. First-order formulations require only first derivatives in the PDE residuals and achieve competitive accuracy, although order reduction alone in the formulation does not consistently improve accuracy under a fixed coefficient budget. These findings motivate the joint consideration of formulation, variable scaling, residual weighting, and feature allocation.

The architectural studies further show that accuracy depends on the random-feature scale, activation, initialization, and space--time basis construction; no single choice dominates across all benchmarks. Adaptive refinement improves several AC results, whereas its benefits are smaller for the tested CH problem. Comparisons with the tested PINN baselines and published results demonstrate substantial accuracy gains, including improvements of one to two orders of magnitude in several reported comparisons.

Future work will address higher-dimensional and long-time problems, adaptive selection of feature scales and field-specific approximation spaces, and residual--error estimates for nonlinear phase-field systems. The present numerical results do not establish the residual--error norm equivalence required for certified accuracy.

%\newpage
\section*{CRediT authorship contribution statement}
CL developed the modeling framework and methodology, implemented the computational code, analyzed the results, and prepared the figures, tables, and original draft of the manuscript.
QW contributed to the conceptualization and methodology development, provided supervision and project administration, secured funding, and participated in the review and editing of the manuscript.
All authors reviewed and approved the final version of the manuscript.

% (Optional) Fill in if needed.

\section*{Declaration of interest}
None.

\section*{Acknowledgements}
The research is supported by NSF award OIA-2242812. QW's research is also partially supported by NSF award DMS-2038080 and DOE award DE-SC0025229. CL gratefully acknowledges Prof.~Jingrun Chen, Prof.~Suchuan Dong, and Xiang Cao for insightful discussions and helpful suggestions regarding the implementation.
The authors gratefully acknowledge the computational resources provided by the Theia high-performance computing cluster at the University of South Carolina, which is supported by NSF award 2320292. We also acknowledge the technical assistance and resources provided by Research Computing at the University of South Carolina (RRID: SCR-027488).
\section*{Data availability}
No data were used in the research described in this article.

\appendix
\begin{appendices} %% 进入附录环境，会自动在页眉显示 Appendices
\counterwithin{table}{section}  % 让表格编号变为 A.1, B.1...
\counterwithin{figure}{section} % 让图片编号变为 A.1, B.1...

% --------------------------------------------------------
\section{Algorithm for gradient-indicator-based adaptive sampling}\label{app:grad_ads}
The gradient-indicator-based adaptive sampling algorithm is presented below.

\begin{algorithm}[htbp]
\caption{Gradient-indicator-based adaptive sampling}
\label{alg:grad-refine}
\begin{algorithmic}[1]
\Require{Subdomains $\{\Omega_s\}$, current sampling sets $\mathcal{P}_s$, trained coefficients $\mathbf{c}_s$, number of extra points $N_{\text{add}}$, quantile $q$, expansion ratio $\alpha$, probe size $N_{\text{probe}}$ (optional)}
\Ensure{Refined sampling sets $\mathcal{P}_s'$}
\State \textbf{Global probe:} Generate a set of candidate points $\mathcal{X}_{\text{probe}}$ uniformly over the entire domain $\Omega = \bigcup_s \Omega_s$.
\State Evaluate the solution $u$ at all probe points using the current DDM solution $\hat{u}(x,t) = \sum_s H_s(x,t) \mathbf{c}_s$.
\State Compute the indicator $I(x) = |\partial_x u|$ (or $\|\nabla u\|$ in 2D) at each probe point.
\State \textbf{Thresholding:} Set $\tau = Q_q(\{I(x)\})$ to the $q$-th quantile of the indicator values.
\State Identify the set of “active” probe points: $\mathcal{A} = \{ x \in \mathcal{X}_{\text{probe}} \mid I(x) \ge \tau \}$.
\State \textbf{Grouping and band formation:}
\If{1D case}
    \State Sort the active points by coordinate and group consecutive points into intervals (clusters).
    \For{each cluster $C$}
        \State Let $x_l = \min C$, $x_r = \max C$.
        \State Expand the interval: $x_l \leftarrow \max(x_{\min}, x_l - \alpha \cdot \text{width}(\Omega))$, $x_r \leftarrow \min(x_{\max}, x_r + \alpha \cdot \text{width}(\Omega))$.
        \State Store the band $B = [x_l, x_r]$.
    \EndFor
\ElsIf{2D case}
    \State Apply connected-component labeling to the active points to obtain clusters $\{C_k\}$.
    \For{each cluster $C$}
        \State Compute its bounding box and dilate it by a factor $\alpha$ (morphological dilation).
        \State Store the resulting rectangular band $B_k$.
    \EndFor
\EndIf
\State \textbf{Local refinement:}
\For{each subdomain $\Omega_s$}
    \State Initialize the list $\mathcal{Q}_s = \emptyset$.
    \For{each global band $B$}
        \State Intersect $B$ with $\Omega_s$: $B_s = B \cap \Omega_s$.
        \If{$B_s$ has positive measure}
            \State Sample $n_s = \max(1, N_{\text{add}} / \text{number of bands})$ new points uniformly within $B_s$.
            \State Append these points to $\mathcal{Q}_s$.
        \EndIf
    \EndFor
    \State $\mathcal{P}_s' \leftarrow \mathcal{P}_s \cup \mathcal{Q}_s$.
\EndFor
\State \Return $\{\mathcal{P}_s'\}$.
\end{algorithmic}
\end{algorithm}

% -----------------------------------------------------------------
 \section{Activation functions: definitions and settings}\label{app:activation} 
For each activation, we provide its mathematical definition and, where applicable, the default hyperparameters used in the PyTorch implementation.

\begin{table}[htbp]
\centering
\caption{Activation functions supported in \texttt{\_phi\_pointwise}.}
\label{tab:activations}
\begin{tabular}{@{}l l c@{}}
\toprule
\textbf{Name (code)} & \textbf{Mathematical expression } \(\phi(z)\) & \textbf{Default parameters} \\
\midrule
\verb|tanh|     & \(\displaystyle \tanh(z) = \frac{e^{z} - e^{-z}}{e^{z} + e^{-z}}\) & — \\
\verb|sigmoid|  & \(\displaystyle \sigma(z) = \frac{1}{1 + e^{-z}}\) & — \\
\verb|sin|      & \(\displaystyle \sin(\pi z)\) & The input is multiplied by \(\pi\). \\
\verb|cos|      & \(\displaystyle \cos(\pi z)\) & The input is multiplied by \(\pi\). \\
\verb|gauss|    & \(\displaystyle e^{-z^{2}}\) & — \\
\verb|exp|      & \(\displaystyle e^{z}\) & — \\
\verb|relu|     & \(\displaystyle \max(0, z)\) & — \\
\verb|relu^2| & \(\displaystyle \max(0, z)^{2}\) & — \\
\verb|leakyrelu^2| & \(\displaystyle \begin{cases} z^{2}, & z \ge 0 \\[4pt] (0.01 \, z)^{2}, & z < 0 \end{cases}\) & Negative slope \(\alpha = 0.01\) (hard-coded). \\
\bottomrule
\end{tabular}
\end{table}

\noindent\textbf{Remarks:}
\begin{itemize}
    \item For \verb|sin| and \verb|cos|, the input is multiplied by \(\pi\) before applying the trigonometric function. This is consistent with the common practice in random feature methods where the frequency is controlled by the random weights, and the \(\pi\) factor ensures a \(2\pi\)-periodic behavior for normalized inputs.
    \item For \verb|leakyrelu^2|, the activation is defined as the square of the Leaky ReLU function. The negative slope is fixed at \(\alpha = 0.01\), which is the default value in PyTorch's \texttt{torch.nn.functional.leaky\_relu}.
    \item The \verb|gauss| activation corresponds to the Gaussian (radial basis function) kernel, computed as \(\exp(-z^2)\) without an additional scaling factor.
\end{itemize}

% ----------------------------------
\section{Results for the Allen--Cahn equation}\label{app:ac1d}

\begin{figure}[htbp]
    \centering
    \includegraphics[width=0.95\linewidth]{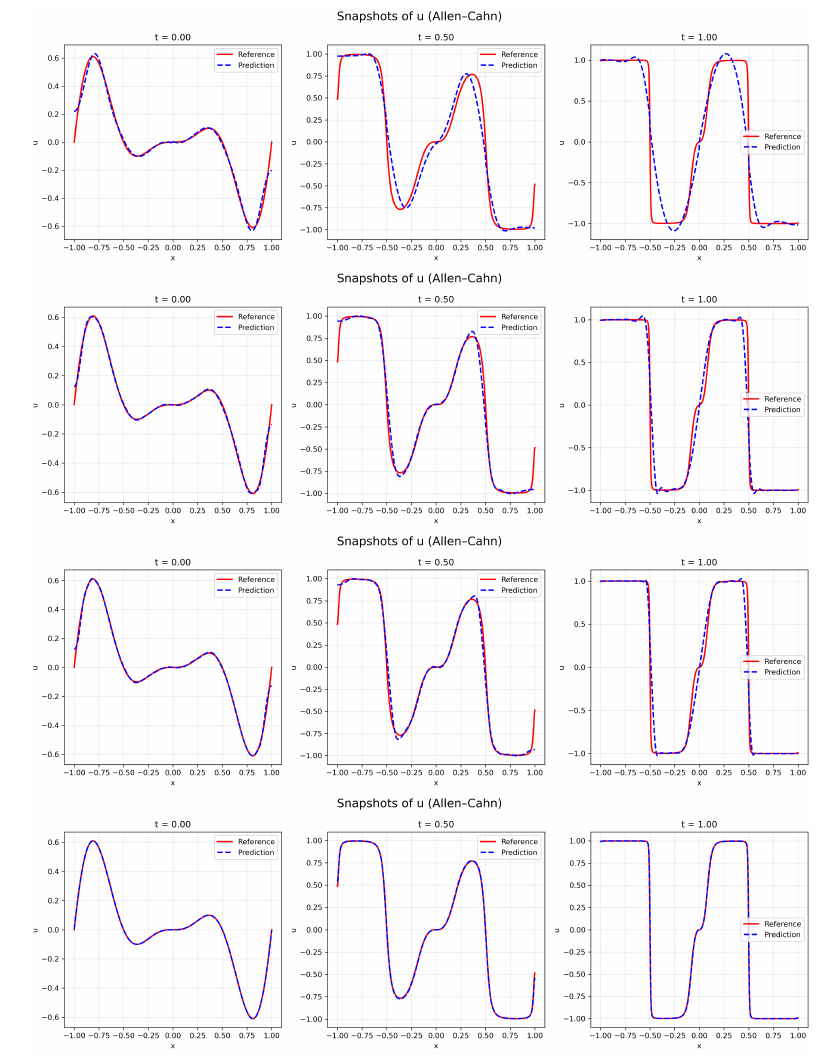}
    \caption{Snapshots at $t=0, 0.5, 1.0$ for the Allen--Cahn equation with $u_0=x^2\sin(2\pi x)$: first row: PINN+Adam, $2\times 200\times 1$, second row: PINN+Adam, $2\times 19\times 19\times 19\times 1$, third row: PINN+Adam, $2\times 128\times 128\times 128\times 1$, fourth row: single-scale RPNN using separable bases, $\tanh$ activation, and normally distributed random features with $R_\xi=1$ and $M=800$.}
    \label{fig:ac1d_rpnn_vs_pinn_sin}
\end{figure}

\begin{figure}[htbp]
    \centering
    \includegraphics[width=0.95\linewidth]{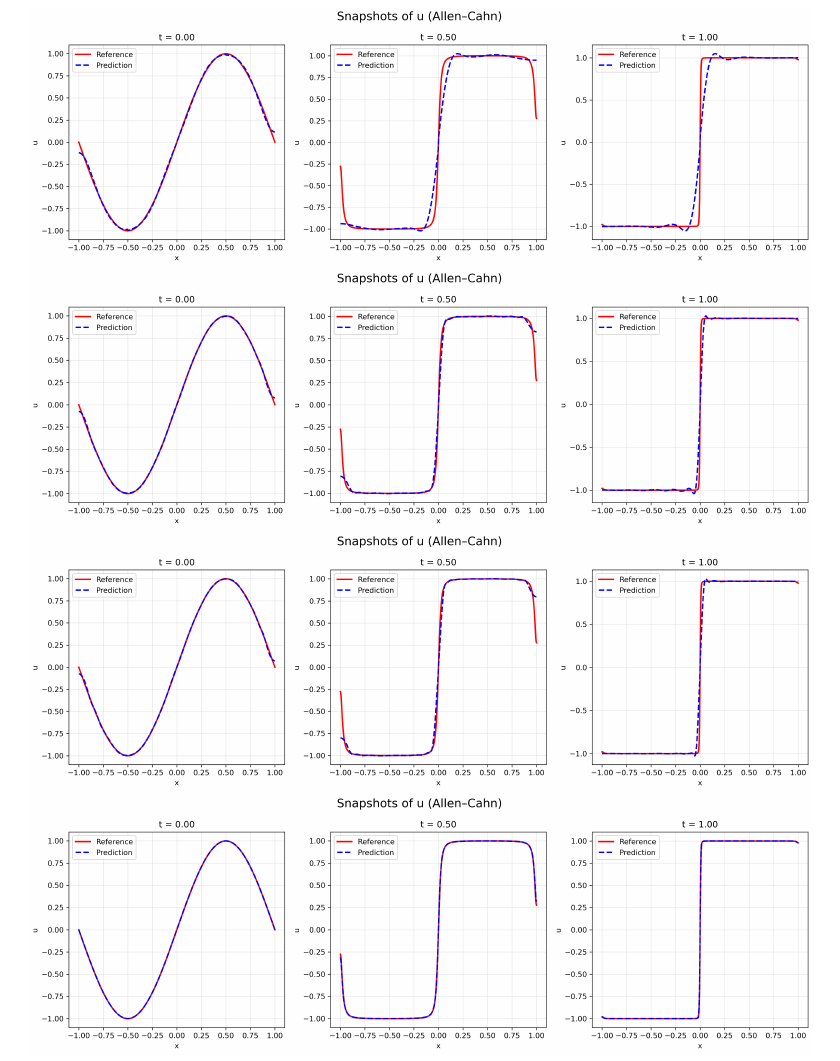}
    \caption{Snapshots at $t=0, 0.5, 1.0$ for the Allen--Cahn equation with $u_0=\sin(\pi x)$: first row: PINN+Adam, $2\times 200\times 1$, second row: PINN+Adam, $2\times 19\times 19\times 19\times 1$, third row: PINN+Adam, $2\times 128\times 128\times 128\times 1$, fourth row: single-scale RPNN using separable bases, $\tanh$ activation, and normally distributed random features with $R_\xi=1$ and $M=800$.}
    \label{fig:ac1d_rpnn_vs_pinn_sinpix}
\end{figure}

The snapshots at $t=0, 0.5, 1.0$ for the Allen--Cahn equation with $u_0=x^2\sin(2\pi x)$ and $u_0=\sin(\pi x)$ are shown in Figure~\ref{fig:ac1d_rpnn_vs_pinn_sin} and Figure~\ref{fig:ac1d_rpnn_vs_pinn_sinpix}.

% ---------------------------------------
\section{Results for the Cahn--Hilliard equation}\label{app:ch1d}
Table~\ref{tab:ch1d_err_rpnn_ablation} shows that residual scaling remains important after variable rescaling. Increasing the reported chemical-potential residual factor from $1$ to $100$ improves every listed RPNN configuration. Variable scaling and equation weighting are distinct design choices.

\begin{figure}[htbp]
    \centering
    \includegraphics[width=0.95\linewidth]{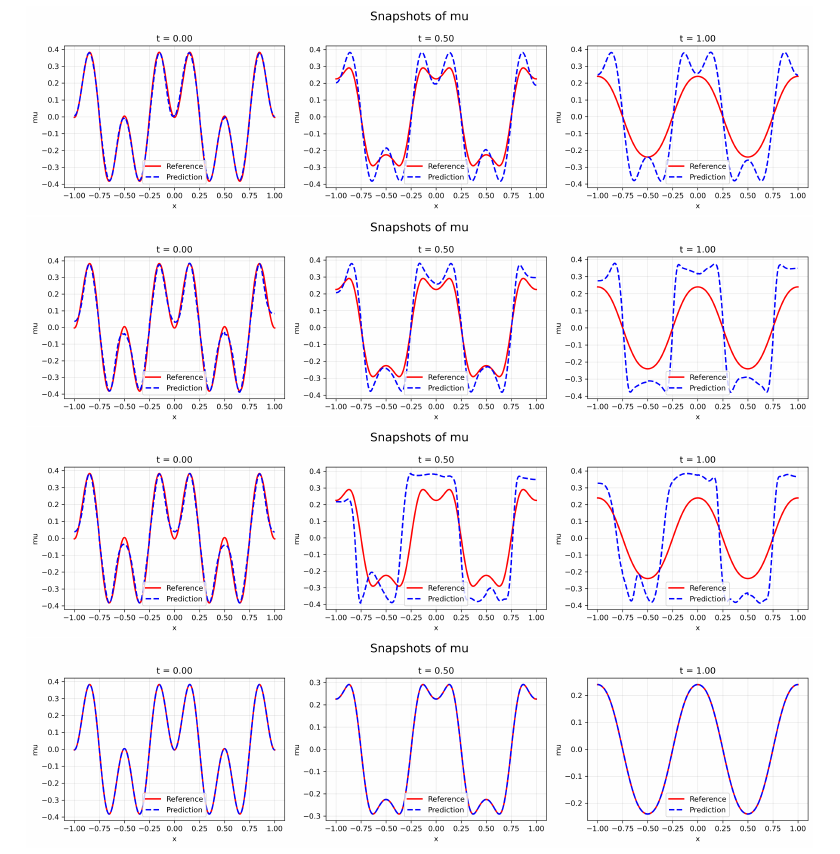}
    \caption{Snapshots of $\mu$ at $t=0, 0.5, 1.0$ for the Cahn--Hilliard equation with $u_0=-\cos(2\pi x)$: first row: PINN-1; second row: PINN-2; third row: PINN-3; fourth row: single-scale RPNN using dense bases, $\tanh$ activation, and uniformly distributed random features with $R_\xi=1$ and $M=800$.}
    \label{fig:ch_rpnn_pinn_cos_uniform_mu}
\end{figure}

Snapshots of $\mu$, as defined in Case 3, for the Cahn--Hilliard equation with $u_0=-\cos(2\pi x)$ are given in Figure~\ref{fig:ch_rpnn_pinn_cos_uniform_mu}. 

% Figure~\ref{fig:3mu_l001} shows the auxiliary variable $\mu$ under the three scalings for $\ell=0.01$: Cases~1 and~2 exhibit a severe imbalance between $u$ and $\mu$, whereas Case~3 displays a well-balanced relation between them.
% \begin{figure}[htbp]
%     \centering
%     \includegraphics[width=0.9\linewidth]{fig/ch1d/mu_3cases/CH_mu_scalings_FDM_001.pdf}
%     \caption{Auxiliary variable $\mu_s$ under the three scalings for $\ell=0.01$. Case 1 (top left) and Case 2 (top right) both show severe imbalance between $u$ and $\mu$; Case 3 (bottom left) shows a well-balanced relation between $u$ and $\mu$; the bottom right panel shows $u$.} %  Physical chemical potential is $s\mu_s$.}
%     \label{fig:3mu_l001}
% \end{figure}

% ---------------- Table B: RPNN w=1 vs w=100 ablation -----------------------------

\begin{table}[htbp]
\centering
\caption{Effect of residual scaling factor \(w_{\mu_s,\text{PDE}}\) on Mixed RPNN structures (Case A: \(u(x,0)=-\cos(2\pi x)\), activation: \(\tanh\)).}
\label{tab:ch1d_err_rpnn_ablation}
\resizebox{\linewidth}{!}{%
\begin{tabular}{lccccc}
\toprule
Model Structure & \(w_{\mu_s, \text{PDE}}\) & \(\text{rel} L^2_u\) & \(L^{\infty}_u\) & \(\text{rel} L^2_{\mu_s}\) & \(L^{\infty}_{\mu_s}\) \\
\midrule
\multirow{2}{*}{Mixed RPNN-d+u} 
& 1   & \(2.936\times10^{-3}\) & \(2.115\times10^{-2}\) & \(7.399\times10^{-3}\) & \(2.729\times10^{-2}\) \\
& 100 & \(\mathbf{8.638\times10^{-4}}\) & \(\mathbf{5.299\times10^{-3}}\) & \(\mathbf{1.024\times10^{-3}}\) & \(\mathbf{1.656\times10^{-3}}\) \\
\midrule
\multirow{2}{*}{Mixed RPNN-d+n} 
& 1   & \(5.528\times10^{-3}\) & \(5.217\times10^{-2}\) & \(1.666\times10^{-2}\) & \(4.665\times10^{-2}\) \\
& 100 & \(\mathbf{4.433\times10^{-3}}\) & \(\mathbf{4.237\times10^{-2}}\) & \(\mathbf{8.925\times10^{-3}}\) & \(\mathbf{1.915\times10^{-2}}\) \\
\midrule
\multirow{2}{*}{Mixed RPNN-s+u} 
& 1   & \(2.107\times10^{-3}\) & \(2.351\times10^{-2}\) & \(8.236\times10^{-3}\) & \(2.602\times10^{-2}\) \\
& 100 & \(\mathbf{8.896\times10^{-4}}\) & \(\mathbf{4.824\times10^{-3}}\) & \(\mathbf{8.633\times10^{-4}}\) & \(\mathbf{1.389\times10^{-3}}\) \\
\midrule
\multirow{2}{*}{Mixed RPNN-s+n}  
& 1 & \(3.055\times 10^{-1}\) & \(1.805\times 10^{0}\) & \(4.051\times 10^{-1}\) & \(8.920\times 10^{-1}\) \\
& 100 & \(\mathbf{7.249\times 10^{-3}}\) & \(\mathbf{5.344\times 10^{-2}}\) & \(\mathbf{5.873\times 10^{-3}}\) & \(\mathbf{1.002\times 10^{-2}}\) \\
\bottomrule
\end{tabular}}
\end{table}

\end{appendices}

\bibliographystyle{plain}
\bibliography{ref}

\end{document}